\documentclass{amsart}
\usepackage[utf8]{inputenc}

\usepackage{amsthm}

\usepackage{mathrsfs}

\usepackage{enumerate}
\usepackage{mathabx}
\usepackage{verbatim}
\usepackage{xcolor}

\usepackage[yyyymmdd,hhmmss]{datetime}

\newcommand{\bigslant}[2]{{\raisebox{.2em}{$#1$}\left/\raisebox{-.2em}{$#2$}\right.}}

\newcommand{\eclass}[2]{[#1]_{\mathscr{#2}}}
\newcommand{\uclass}[1]{\eclass{#1}{U}}

\newtheorem{theorem}{Theorem}[section]
\newtheorem{observation}[theorem]{Observation}
\newtheorem{claim}{Claim}[theorem]
\newtheorem{lemma}[theorem]{Lemma}
\newtheorem{corollary}[theorem]{Corollary}
\theoremstyle{definition}
\newtheorem{definition}[theorem]{Definition}
\newtheorem{notion}[theorem]{Notion}
\theoremstyle{remark}

\newcommand{\Symm}{\operatorname{Sym}}

\title{On proper colorings of functions}
\author{Tamás Csernák}
\address{HUN-REN Alfréd Rényi Institute of Mathematics, Reáltanoda u 13--15 Budapest, 1053 Hungary}
\email{tamas@csernak.com}


\begin{document}

\begin{abstract} 
We investigate the infinite version of the $k$-switch problem of 
Greenwell and Lovász.

Given infinite cardinals ${\kappa}$ and ${\lambda}$, 
for functions  $x,y\in {}^{\lambda}\kappa $ we say that they are \emph{totally different} 
if $x(i)\ne y(i)$ for each $i\in {\lambda}$.
A function $F:{}^{\lambda}\kappa \longrightarrow {\kappa} $ is a 
\emph{proper coloring} if 
$F(x)\ne F(y)$ whenever  $x$ and $y$ are totally different elements of 
 ${}^\lambda{\kappa} $. 
We say that 
\begin{enumerate}[(i)]
\item $F$ is {\em weakly uniform} iff there are pairwise totally different
functions $\{r_{\alpha}:{\alpha}<{\kappa}\}\subset {}^{\lambda}{\kappa}$
such that $F(r_{\alpha})={\alpha}$;
\item $F$ is \emph{tight} if 
there is no proper coloring $G:{}^{\lambda}\kappa \longrightarrow {\kappa}$
such that there is exactly one $x\in {}^{\lambda}{\kappa}$ 
with  $G(x)\ne F(x)$.
\end{enumerate}

We show that  given  a  proper coloring $F:{}^{\lambda}{\kappa}\to {\kappa}$,
the following statements are equivalent 
\begin{enumerate}[(i)]
 \item $F$ is weakly uniform,
 \item there is a ${\kappa} ^{+}$-complete ultrafilter $\mathscr{U}$  on 
 ${\lambda}$ and there is a permutation ${\pi}\in \Symm({\kappa})$
 such that for each $x\in {}^{\lambda}{\kappa}$ we have 
 \begin{displaymath}
 F(x)={\pi}({\alpha})\ \Longleftrightarrow \ \{i\in {\lambda}: x(i)={\alpha}\}
 \in \mathscr{U}.
 \end{displaymath}

\end{enumerate}
We also show that there are tight proper colorings    
which cannot be obtained such a way. 

\end{abstract}

\maketitle

\section{Introduction}

The well-known $k$-switch theorem of Greenwell and Lov\'asz \cite{GLL} has stated the following: 
{\em Assume that there is a light with $n\ge 3$ states, which is connected to finitely many switches, each having $n$ positions, in such a way that    
if the positions of all switches are changed, then the state of the light also changes. Then there is a switch  such that 
 the position of that switch determines the state of the light.}


More formally, we denote the possible configurations of the switches by elements of ${}^{k}{n}$, and let $F:{}^{k}{n}\longrightarrow n$ be a function such that if $x,y\in {}^{k}{n}$ with $x(i)\neq y(i)$ for all $i\in k$ then  $F(x)\neq F(y)$. In this case, there is a $j\in k$ and a permutation $\pi \in Sym({n})$, such that $F(x)=\pi (x(j))$ for all $x\in {}^{k}{n}$.

In this paper we try to generalize this theorem to the infinite case, involving functions $F:{}^{\lambda}\kappa \longrightarrow \kappa $, or more general 
$F:{}^{\lambda}\kappa \longrightarrow \mu $, where $\kappa ,\lambda ,\mu $ are cardinals. 
For any  $x,y\in {}^{\lambda}\kappa $ let $$\Delta (x,y)=\{ i\in \lambda:x(i)\neq y(i)\}.$$
For $x,y\in {}^{\lambda}\kappa $ we say that they are \emph{totally different}, and denote it by $x\perp y$, if 
$\Delta(x,y)={\lambda}$.
For a set $A\subseteq \lambda $, we say $x$ and $y$ are \emph{totally different on $A$} and denote it by $x\perp _{A}y$, if $A\subset \Delta(x,y)$, and they \emph{coincide on $A$}, denoted by $x\equiv _{A}y$ if 
$\Delta(x,y)\cap A=\emptyset$.
  So $x\perp _{\Delta (x,y)}y$ and $x\equiv _{\lambda -\Delta (x,y)}y$.

\begin{definition}
    A function $F:{}^{\lambda}\kappa \longrightarrow \mu $ is a \emph{proper coloring} if $F(x)\neq F(y)$ for each $x,y\in {}^{\kappa}\lambda $ with $x\perp y$. 

    A proper coloring $G:{}^{\lambda}\kappa \longrightarrow \mu $ is \emph{trivial} if ${\mu}={\kappa}$ and there is an $i\in {\lambda}$ and ${\pi}\in Sym({\kappa})$, such that   $G(x)={\pi}(x(i))$ for each $x\in ^{\lambda}{\kappa}$.   
  \end{definition}

Komjáth and  Totik  \cite{KoTo} 
considered the case when there are infinitely many switches and 
every switch has $n$ states for some $3\le n<{\omega}$,
and they obtained the following generalization of the original result of 
Greenwell and Lovász.

\begin{theorem}[{\cite[Section 3]{KoTo}}] \label{PRO_TmKoTo}
    Assume that ${\lambda}$ is a (finite or infinite) cardinal and   $n\ge 3$ is a natural number. 
 A  function $F:{}^{\lambda }{n}\longrightarrow n$ is a proper coloring if and only if there is an  ultrafilter $\mathscr{U}$ on $\lambda $ and 
 there is a permutation $\pi \in Sym({n})$, such that for all $k<n$ and $x\in {}^{\lambda}n$,
 \begin{displaymath}
  \tag{$\dag$}F(x)=\pi (k)\ \Longleftrightarrow\ \{ i\in \lambda\ |\ x(i)=k\} \in \mathscr{U}.
 \end{displaymath} 
\end{theorem}

Our aim is to generalize this result for proper colorings 
$F:{}^{\lambda}\kappa \longrightarrow \kappa $ or $F:{}^{\lambda}\kappa \longrightarrow \mu $, where $\kappa $ is infinite.



\begin{notion}
Given any set $A$, we denote by $Sym(A)$ the set of all bijections $\pi :A\longrightarrow A$.
\end{notion}

\begin{notion}
For a set $A$, an ultrafilter $\mathscr{U}\subset \mathcal{P}(A)$ and a cardinal $\kappa $, we say that $\mathscr{U}$ is $\kappa $-complete is for any $\mathscr{F}\subseteq \mathscr{U}$, $|\mathscr{F}|< \kappa $, we have $\bigcap{\mathscr{F}}\in \mathscr{U}$. Hence the ultrafilter is $\kappa ^{+}$-complete, if it is also closed for intersections of $\kappa $ many sets.
\end{notion}

\section{Unrestricted colorings}

Let $F:{}^{\lambda}\kappa \longrightarrow \mu $ be a proper coloring, where $\kappa ,\mu $ are infinite cardinals. 
Furthermore, for $\alpha \in \kappa $ we denote by $c_{\alpha }$ the constant ${\alpha}$ function from ${}^{\lambda}{\kappa}$, i.e.  $c_{\alpha }(i)=\alpha $ for all $i\in \lambda$. 
For all $\alpha \ne\beta \in \kappa $, we have $c_{\alpha }\perp c_{\beta}$, so $F(c_{\alpha })\neq F(c_{\beta })$. The images of constant functions must be all different, so  ${\kappa}\le |ran(F)|\le {\mu}$.
i.e,    in the case $\mu <\kappa $ there are no proper colorings.

In the case $\mu =\kappa $, we have proper colorings in the form 
$F(x)={\pi}(x(i))$ for any fixed $i\in \lambda $ and ${\pi}\in Sym({\kappa})$.
 These colorings are considered as the \emph{trivial} ones, but in general we can ask if every 
proper coloring $F:{}^{\lambda}\kappa \longrightarrow \kappa $  can be given in a "nice" form, just like in the finite range case. The next easy observation says that in general there is no  "nice" form for it, as we can extend any partial proper coloring in some way.

\begin{observation}
 Assume that  $\mathscr{A}\subseteq {}^{\lambda }{\kappa } $ and $G:\mathscr{A}\longrightarrow \kappa$ is a {\em partial proper coloring} i.e. if $x,y\in \mathscr{A}$ with $x\perp y$, then $G(x)\neq G(y)$.  Then there is a proper coloring $F:{}^{\lambda}\kappa \longrightarrow \kappa $, such that for each $x,y\in \mathscr{A}$, $F(x)=F(y)\Leftrightarrow G(x)=G(y)$.
\end{observation}

\begin{proof}
 Let $A,B\subseteq \kappa $ with $A\cap B=\emptyset $ and $|A|=|B|=\kappa $. Let $h_{1}:\kappa \longrightarrow A$ and $h_{2}:\kappa \longrightarrow B$ be bijections. We define $F$ in the following way: If $x\in \mathscr{A}$, then $F(x)=h_{1}(G(x))$, and if $x\not\in \mathscr{A}$, then $F(x)=h_{2}(x(0))$. 

 The coloring $F$  clearly satisfies the requirements.
 \end{proof}


To make proper colorings more manageable, we should impose additional properties. For instance, in the next section, we will constrain their values on a relatively small set, the set of constant functions.

\section{Uniform colorings}

\begin{definition}\label{PRO_Def1}
 A proper coloring $F: {}^{\lambda}\kappa \longrightarrow \kappa $  is \emph{strongly uniform} if $F(c_{\alpha })=\alpha $ for all $\alpha \in \kappa$.
\end{definition}

\begin{lemma}\label{PRO_Lm1}
 If $F:{}^{\lambda}\kappa \longrightarrow \kappa $ is a strongly uniform proper coloring, $S\subseteq \kappa $ and $x\in {}^{\lambda}\kappa $ is such that $x(i)\in S$ for all $i\in \lambda $, then $F(x)\in S$.
\end{lemma}

\begin{proof}
 Choose any $\beta \in \kappa $, such that $\beta \not\in S$. Then $x\perp c_{\beta }$, so $F(x)\neq F(c_{\beta })=\beta $. Hence, $F(x)$ differs from all $\beta \not\in S$, thus $F(x)\in S$.
\end{proof}

Now we focus on the case where $S$ is a finite subset of $\kappa $
using the theorem of Komjáth and Totik.  

\begin{lemma}\label{PRO_Lm2}
 If $F:{}^{\lambda}\kappa \longrightarrow \kappa $ is a strongly uniform proper coloring, then there is an ultrafilter $\mathscr{U}$ on $\lambda $, such that for every $x\in {}^{\lambda }{\kappa }$ with $|Ran(x)|<\omega $ and $\alpha \in \kappa$, we have $$F(x)=\alpha\ \Longleftrightarrow\ \{ i\in \lambda\ |\ x(i)=\alpha \} \in \mathscr{U}.$$
\end{lemma}

\begin{proof}
 For any $S\subseteq \lambda$, $|S|<\omega $ let $F_{S}=F|_{^{\lambda }S}$. 
 By Lemma \ref{PRO_Lm1} we know that $F_{S}$ maps into $S$, so $F_{S}:{}^{\lambda }S\longrightarrow S$ is a proper coloring on $^{\lambda }S$. For every $|S|\ge 3$ 
 by Theorem \ref{PRO_TmKoTo}, using the fact that $F_S(c_{\alpha})={\alpha}$
 for each ${\alpha}\in S$,  there is an  ultrafilter $\mathscr{U}_{S}$  on ${\lambda}$
 such that  
 \begin{displaymath}
  \tag{$\dag_S$}\forall x\in {}^{\lambda}S\ \big(\ F(x)=k\ \Longleftrightarrow\ \{ i\in \lambda\ |\ x(i)=k\} \in \mathscr{U}_S\big).
\end{displaymath}

 Let  $S,T$ be  finite subsets of $\kappa $ with at least 3 elements and $S\subseteq T$. Pick $a\ne b \in S$.
Let $A\in \mathscr{U}_T$. Define $x\in {}^{\lambda}S$ such that 
$x\equiv_A c_a$ and $x\equiv_{{\lambda}\setminus A}c_b$.
Then $F_T(x)=a$ by $(\dag_T)$. Thus, $A\in \mathscr{U}_S$ by $(\dag_S)$.
Hence  $\mathscr{U}_{T}\subseteq \mathscr{U}_{S}$, but both are ultrafilters on $\lambda $, so $\mathscr{U}_{T}=\mathscr{U}_{S}$.

In general, if $S,T$ are any finite subsets of $\kappa $ with at least 3 elements, then $\mathscr{U}_{S}=\mathscr{U}_{S\cup T}=\mathscr{U}_{T}$. So there is a unique ultrafilter $\mathscr{U}$ on $\lambda $, such that $\mathscr{U}_{S}=\mathscr{U}$ for all at least three element finite subsets $S$ of ${\kappa}$. 

 
 Now consider an arbitrary  $x\in {}^{\lambda }{\kappa }$ with $|Ran(x)|<\omega $. 
 Let $S$ be a finite subset of ${\kappa}$ such that 
 $|S|\ge 3$ and $Ran(x)\subset S$.
 
 Then $F(x)=F_S(x)$ and $F_S(x)=\alpha$ iff 
 $\{i<{\lambda}:x(i)={\alpha}\}\in \mathscr{U}_S=\mathscr{U}$.
\end{proof}

\begin{definition}\label{PRO_Def2}
 For any strongly uniform proper coloring $F:{}^{\lambda}\kappa \longrightarrow \kappa $ we define the \emph{ultrafilter corresponds to $F$} as the unique $\mathscr{U}$ ultrafilter on $\lambda $ such that for every $x\in {}^{\lambda }{\kappa }$ with $|Ran(x)|<\omega $ and $\alpha \in \kappa$, we have $$F(x)=\alpha \ \Longleftrightarrow\  \{ i\in \lambda\ |\ x(i)=\alpha \} \in \mathscr{U}.$$ By Lemma  \ref{PRO_Lm2} such an   ultrafilter exists and is clearly unique.
\end{definition}

\begin{theorem}\label{PRO_Thm3}
 If $F:{}^{\lambda}\kappa \longrightarrow \kappa $ is a  strongly uniform proper coloring and $\mathscr{U}$ is the corresponding ultrafilter, then for all $x\in {}^{\lambda }{\kappa } $ and $\alpha \in \kappa $ we have $$F(x)=\alpha \ \Longleftrightarrow\ \{ i\in \lambda\ |\ x(i)=\alpha \} \in \mathscr{U}$$ (so we can omit the condition $|Ran(x)|<\omega $ from Lemma  \ref{PRO_Lm2})
\end{theorem}

\begin{proof}
 First suppose $\{ i\in \lambda |x(i)=\alpha \} \not\in \mathscr{U}$. Choose any $\beta \in \kappa $ with $\beta \neq \alpha $, and let $y\in {}^{\lambda }{\kappa } $ be such that $y(i)=\alpha $ if $x(i)\neq \alpha $ and $y(i)=\beta $ if $x(i)=\alpha $. Then $Ran(y)=\{ \alpha ,\beta \}$ is finite, and $\{ i\in \lambda |y(i)=\alpha \} =\{ i\in \lambda |x(i)\neq \alpha \} \in \mathscr{U}$, so $F(y)=\alpha $. On the other hand we constructed $y$ such that $y\perp x$, and so $F(x)\neq F(y)=\alpha $ in this case.
 
 Now suppose $\{ i\in \lambda |x(i)=\alpha \} \in \mathscr{U}$. Then for any $\beta \neq \alpha $ we have $\{ i\in \lambda |x(i)=\beta \} \not\in \mathscr{U}$, so as we proved in the previous paragraph $F(x)\neq \beta $, thus $F(x)=\alpha $ is the only possible value.
\end{proof}

\begin{lemma}\label{PRO_Lm3}
 If $F:{}^{\lambda}\kappa \longrightarrow \kappa $ is a strongly uniform proper coloring then the corresponding ultrafilter must be $\kappa ^{+}$-complete.
\end{lemma}

\begin{proof}
 Suppose, for contradiction, that $\mathscr{U}$ is not $\kappa ^{+}$-complete, and let $A_{\alpha }\in \mathscr{U}$ for $\alpha \in \kappa $ be sets with $\bigcap_{\alpha \in \kappa }A_{\alpha } \not\in \mathscr{U}$. Now define $x\in {}^{\lambda }{\kappa } $ in the following way: Let $x(i)=min\{ \alpha |i\not\in A_{\alpha }\}$ if $i\not\in \bigcap_{\alpha \in \kappa }A_{\alpha }$ and $x(i)=0$ if $i\in \bigcap_{\alpha \in \kappa }A_{\alpha }$. Also let $B_{\alpha }=\{ i\in \lambda |x(i)=\alpha \}$. If $\alpha \neq 0$, then by definition of $x$, we have $A_{\alpha }\cap B_{\alpha }=\emptyset $, and since $A_{\alpha }\in \mathscr{U}$, $B_{\alpha }\not\in \mathscr{U}$. In the case ${\alpha}=0$, we have $B_{0}=(\lambda -A_{0})\cup \bigcap_{\alpha \in \kappa }A_{\alpha }$ and neither $\lambda -A_{0}$ nor $\bigcap_{\alpha \in \kappa }A_{\alpha }$ are in $\mathscr{U}$, so $B_{0}\not\in \mathscr{U}$. Then for any $\alpha \in \kappa $, since $B_{\alpha }\not\in \mathscr{U}$, by Theorem \ref{PRO_Thm3} $F(x)\neq \alpha $, so $F(x)$ can not have any value and that is a contradiction.
\end{proof}

This lemma gives a strict restriction concerning   $\mathscr{U}$ and $F$.
 Clearly all principal ultrafilters are $\kappa ^{+}$-complete, 
 but non-principal  $\kappa ^{+}$-complete ultrafilters exist 
 on ${\lambda}$
 only if there is a measurable cardinal between $\kappa $ and $\lambda $. 
 Measurable cardinals are very large, so in many cases, 
 the following corollary can be useful.

\begin{corollary}\label{PRO_Cor1}
If $\lambda $ is below the smallest measurable cardinal or $\lambda \le \kappa $, and $F:{}^{\lambda}\kappa \longrightarrow \kappa $ is a strongly uniform proper coloring, then there is an $i\in \lambda $, such that $F(x)=x(i)$ for all $x\in {}^{\lambda }{\kappa }$. 
\end{corollary}

\begin{proof}
 By Lemma  \ref{PRO_Lm3}, the corresponding ultrafilter $\mathscr{U}$ must be $\kappa ^{+}$-complete, and thus $\omega_{1}$-complete on $\lambda $. If $\lambda $ is below the smallest measurable cardinal, $\mathscr{U}$ must be principal, and in case $\lambda \le \kappa $ it is clearly principal, so there is $i\in \lambda $ such  that for all $X\subseteq \lambda $, $X\in \mathscr{U}\Leftrightarrow i\in X$. Then by Theorem \ref{PRO_Thm3}, for any $x\in {}^{\lambda }{\kappa } $ and $\alpha \in \kappa $ we have $F(x)=\alpha \Leftrightarrow  \{ j\in \lambda |x(j)=\alpha \} \in \mathscr{U} \Leftrightarrow i\in \{ j\in \lambda |x(j)=\alpha \} \Leftrightarrow x(i)=\alpha $, thus $F(x)=x(i)$.
\end{proof}

Now we have seen that strongly uniform proper colorings can be written in a nice form. But if we remember the definition with  switches, we can see that strong uniformity is not a natural assumption in this case, since there is no canonic bijection between the states of the switches and the colors of the lamp. To overcome this, we define a more general property  that can be reverted to the strongly uniform case.

\begin{definition}\label{PRO_Def3}
 A proper coloring $F:{}^{\lambda}\kappa \longrightarrow \kappa $ is \emph{weakly uniform} if there is a sequence $(r_{\alpha })_{\alpha <\kappa }$ in ${}^{\lambda }{\kappa } $ such that 
\begin{enumerate}[1.]
\item   $F(r_{\alpha })=\alpha $ for $\alpha \in \kappa $,
\item   $r_{\alpha }\perp r_{\beta }$ for $\alpha \neq \beta \in \kappa $.
\end{enumerate} 
(In other words, the graph $(\mathbf K_{\kappa})^{\lambda}$ contains 
a clique $\{r_{\alpha}:{\alpha}<{\kappa}\}$ such that $F(r_{\alpha})={\alpha}$.)
\end{definition}

It is clear that strongly uniform proper colorings are weakly uniform using $r_{\alpha }=c_{\alpha }$. 


\begin{observation}
If $n$ is a natural number, then every proper coloring
$F:{}^{\lambda}n\to n$ is weakly uniform. 
\end{observation}
Indeed, $|\{F(c_i):i<n\}|=n$ because $c_i\perp c_j$
for $i<j<n$, and so $\{F(c_i):i<n\}=\{0,1,\dots,n-1\}$.

\medskip

Now we try to characterize weakly uniform proper colorings . 


\begin{theorem}\label{PRO_Thm4}
A function $F:{}^{\lambda}\kappa \longrightarrow \kappa $ is a weakly uniform proper coloring  if and only if  there is a $\kappa ^{+}$-complete ultrafilter  $\mathscr{U}$ on $\lambda $ and a permutation $\pi \in Sym(\kappa )$ such that for all $x\in {}^{\lambda }{\kappa } $ and $\alpha \in \kappa $ we have $$F(x)=\pi (\alpha )\Leftrightarrow \{i\in \lambda |x(i)=\alpha \} \in \mathscr{U}.$$
\end{theorem}

\begin{proof}
Assume first that $\mathscr{U}$ is a  $\kappa ^{+}$-complete ultrafilter   on $\lambda $ and $\pi \in Sym(\kappa )$ is a permutation.

Since $\mathscr U$ is ${\kappa}^{+}$-complete,  for each $x\in {}^{\lambda}{\kappa}$ there is 
a color ${\alpha}$ such that $\{i\in \lambda |x(i)=\alpha \} \in \mathscr{U}$.
Thus, $F(x)={\pi}({\alpha})$ is defined. 

If $x\perp y$, $F(x)={\pi}({\alpha})$ and 
$F(y)={\pi}({\beta})$, then 
$\{i: x(i)={\alpha}\land y(i)={\beta}\}\in \mathscr{U}$
and there is no $i$ with $x(i)=y(i)$, so 
${\alpha}\ne {\beta}$ and so $F(x)\ne F(y)$. Thus $F$
is a proper coloring. 

Finally, $F(c_{\pi^{-1}({\alpha})})={\pi}(\pi^{-1}({\alpha}))={\alpha}$,
so $F$ is weakly uniform as well. 

\medskip

To prove the other direction 
  assume that the sequence  $(r_{\alpha })_{\alpha <\kappa }\subset {}^{\lambda}{\kappa}$ witnesses that 
$F$ is weakly uniform. 

For each $i<{\lambda}$ define the function $h_i:{\kappa}\to {\kappa} $
by the formula 
\begin{displaymath}
  h_i({\alpha})=r_{\alpha}(i). 
\end{displaymath}

Since $r_{\alpha}\perp r_{\beta}$ for ${\alpha}<{\beta}<{\kappa}$,
the functions $h_i$ are injective. 

  Let $H:{}^{\lambda}\kappa \longrightarrow {}^{\lambda }{\kappa }$ be the function defined as $(H(x))(i)=h_{i}(x(i))$. Since all the functions $h_{i}$ are one-to-one, the function  $H$ is also one-to-one, preserves the totally different relation and its range 
  is the set $ran(H)=\bigtimes _{i\in \lambda }Ran(h_{i})$.
 
 Now let $G=F\circ H$. Then $G:{}^{\lambda}\kappa \longrightarrow \kappa $ is a proper coloring. Indeed, let $x,y\in  {}^{\lambda }{\kappa }$ with $x\perp y$. Then $H(x)\perp H(y)$, so $G(x)=F(H(x))\neq F(H(y))=G(y)$. Now look at the values of $H$ on the constant functions. Let $\alpha \in \kappa $ and $i\in \lambda $. Then $(H(c_{\alpha}))(i)=h_{i}(c_{\alpha }(i))=h_{i}(\alpha )=r_{\alpha }(i)$, thus $H(c_{\alpha })=r_{\alpha }$. Then $G(c_{\alpha })=F(H(c_{\alpha }))=F(r_{\alpha })=\alpha $, so $G$ is strongly uniform.
 
 Let $\mathscr{U}$ be the ultrafilter corresponding to $G$. By Lemma  \ref{PRO_Lm3} $\mathscr{U}$ is $\kappa ^{+}$-complete.  
 \begin{claim}
  For any $x\in {}^{\lambda }{\kappa }$ and $\alpha \in \kappa $ we have $$F(x)=\alpha \Leftrightarrow \{ i\in \lambda |x(i)=r_{\alpha }(i)\} \in \mathscr{U}.$$ 
 \end{claim}
 
 \begin{proof}[Proof of the Claim]
 First we show this for $x\in Ran(H)$. Then 
 \begin{displaymath}
     F(x)=F(H(H^{-1}(x)))=G(H^{-1}(x)).
 \end{displaymath}
 By Theorem \ref{PRO_Thm3}, we have 
 \begin{displaymath}
    F(x)=\alpha \Leftrightarrow G(H^{-1}(x))=\alpha \Leftrightarrow \{i\in \lambda |(H^{-1}(x))(i)=\alpha \} \in \mathscr{U}.
 \end{displaymath} 
 On the other hand, clearly 
 \begin{displaymath}
   H^{-1}(x)(i)=\alpha \Leftrightarrow h_{i}^{-1}(x(i))=\alpha \Leftrightarrow x(i)=h_{i}(\alpha )\Leftrightarrow  x(i)=r_{\alpha }(i). 
  \end{displaymath} 
   For all $i\in \lambda$, so
   \begin{displaymath}
      \{i\in \lambda |(H^{-1}(x))(i)=\alpha \} =\{ i\in \lambda |x(i)=r_{\alpha }(i)\}. 
   \end{displaymath}

 So the claim holds for $x\in Ran(H)$.

 Now let $x\in {}^{\lambda }{\kappa }$ be arbitrary and $\alpha \in \kappa $. First suppose $\{ i\in \lambda |x(i)=r_{\alpha }(i)\} \not\in \mathscr{U}$. Construct a $y\in {}^{\lambda }{\kappa } $ in the following way: Let $y(i)=r_{\alpha }(i)$ if $x(i)\neq r_{\alpha }(i)$. If $x(i)=r_{\alpha }(i)$ then choose $y(i)\in Ran(h_{i})$ with $y(i)\neq x(i)$ (we can choose such $y(i)$ since the range of a one-to-one mapping from an infinite set must be infinite). Then $y\in Ran(H)=\bigtimes _{i\in \lambda }Ran(h_{i})$ and $\{ i\in \lambda |y(i)=r_{\alpha }(i)\} = \{ i\in \lambda |x(i)\neq r_{\alpha }(i)\} \in \mathscr{U}$, so $F(y)=\alpha $. On the other hand, by the choice of $y$, we have $y\perp x$, so $F(x)\neq F(y)=\alpha $. Now suppose $\{ i\in \lambda |x(i)=r_{\alpha }(i)\} \in \mathscr{U}$. Then for all $\beta \in \kappa $, $\beta \neq \alpha $, we have $\{ i\in \lambda |x(i)=r_{\beta }(i)\} \cap \{ i\in \lambda |x(i)=r_{\alpha }(i)\} =\emptyset $, since $r_{\alpha }\perp r_{\beta }$. Hence $\{ i\in \lambda |x(i)=r_{\beta }(i)\} \not\in \mathscr{U}$, so $F(x)\neq \beta $, thus $F(x)=\alpha $ is the only possible value.
Thus, we proved the Claim. 
\end{proof}

 Let $\pi (\beta )=F(c_{\beta })$ for $\beta \in \kappa $. Since the constant functions are pairwise totally different and $F$ is a proper coloring, $\pi $ is clearly one-to-one. We need to show that $\pi $ is onto. Choose any $\alpha \in \kappa $. The value of $r_{\alpha }$ defines a partition of $\lambda $ into $\kappa $ pieces, and $\mathscr{U}$ is $\kappa ^{+}$ complete, so for some $\beta \in \kappa $ we have $\{ i\in \lambda |r_{\alpha }(i)=\beta \} \in \mathscr{U}$. Since $\{ i\in \lambda |c_{\beta }(i)=r_{\alpha }(i)\} \in \mathscr{U}$, we have $\pi (\beta )=F(c_{\beta })=\alpha $, so $\pi $ is onto. On the other hand, we have for all $\beta \in \kappa$ that $\{ i\in \lambda |r_{\pi (\beta )}(i)=\beta \} \in \mathscr{U}$, so for any $x\in {}^{\lambda }{\kappa }$ and $\beta \in \kappa$ we have $$F(x)=\pi (\beta ) \Leftrightarrow \{ i\in \lambda |x(i)=r_{\pi (\beta ) }(i)\} \in \mathscr{U} \Leftrightarrow \{ i\in \lambda |x(i)=\beta \} \in \mathscr{U}.$$
 
Thus, we completed the proof of the theorem. 
\end{proof}

Theorem \ref{PRO_Thm4} clearly yields the following corollaries.

\begin{corollary}\label{PRO_Cor2}
 If $\lambda $ is below the smallest measurable cardinal or $\lambda \le \kappa$, and $F:{}^{\lambda}\kappa \longrightarrow \kappa $ is weakly uniform a proper coloring, 
 then there are $i\in \lambda $ and $\pi \in Sym(\kappa )$, such that $F(x)=\pi (x(i))$ for all $x\in {}^{\lambda}{\kappa }$. 
\end{corollary}

\begin{proof}
 Let $\mathscr{U}$ be the $\kappa ^{+}$-complete ultrafilter defined by Theorem \ref{PRO_Thm4}, and that is also $\omega _{1}$-complete on $\lambda $. Since $\lambda $ is below the smallest measurable cardinal or $\lambda\le \kappa $, $\mathscr{U}$ must be principal, so there is $i\in \lambda $ that for all $X\subseteq \lambda $, $X\in \mathscr{U}\Leftrightarrow i\in X$. Then by Theorem \ref{PRO_Thm4}, for any $x\in {}^{\lambda }{\kappa } $ and $\alpha \in \kappa $ we have $F(x)=\pi (\alpha ) \Leftrightarrow  \{ j\in \lambda |x(j)=\alpha \} \in \mathscr{U} \Leftrightarrow i\in \{ j\in \lambda |x(j)=\alpha \} \Leftrightarrow x(i)=\alpha $, thus $F(x)=\pi (x(i))$ by the permutation $\pi \in Sym(\kappa )$ from Theorem \ref{PRO_Thm4}.
\end{proof}

\section{Minimal and tight colorings}

In this section, we try to give properties witch are weaker than strong or weak uniformity, but  are strong enough  to ensure that some proper colorings have a nice form. We also look at the more general case when $F:{}^{\lambda}\kappa \longrightarrow \mu $ for some ${\mu}\ge {\kappa}$.

\begin{definition}\label{PRO_Def4}
 For proper colorings $F,G:{}^{\lambda}\kappa \longrightarrow \mu $  we write that $G\le F$ if $G(x)\le F(x)$ for all $x\in {}^{\lambda }{\kappa } $, and $G<F$ if $G\le F$ and $G\neq F$. The relation $\le$ is a partial order on proper colorings. A proper coloring  $F:{}^{\lambda}\kappa \longrightarrow \mu $ is \emph{minimal} if there is no proper coloring $G:{}^{\lambda}\kappa \longrightarrow \mu $ with $G<F$.
\end{definition}

Although the relation $\le $ is not well-founded we can obtain the following result. 
\begin{lemma}\label{PRO_Lm4}
 If $F:{}^{\lambda}\kappa \longrightarrow \mu $ is a proper coloring, then there is a minimal proper coloring $G:{}^{\lambda}\kappa \longrightarrow \mu $ with $G\le F$.
\end{lemma}

\begin{proof}
 Suppose for contradiction that there is no such $G$. Let 
 \begin{displaymath}
  L_{F}=\{ G\ |\ G:{}^{\lambda}\kappa \longrightarrow \mu 
  \text{ is a proper coloring, }G\le F\}.
 \end{displaymath}
   If it does not contain any minimal proper coloring then let $H:L_{F}\longrightarrow L_{F}$ be a function such that $H(G)<G$ for all $G\in L_{F}$. Now we construct a strictly decreasing transfinite sequence $\langle G_{\alpha }:\alpha\in On\rangle$ in $L_{F}$. Let $G_{0}=F$ and $G_{\alpha +1}=H(G_{\alpha })$ for $\alpha \in On $. If $\alpha $ is a limit ordinal then let $G_{\alpha }$ be the function with $G_{\alpha }(x)=min\{ G_{\beta }(x):\beta <\alpha \}$ for all $x\in {}^{\lambda }{\kappa } $. We need to see that this is a proper coloring. Let $x,y\in {}^{\lambda }{\kappa } $ with $x\perp y$. By the definition of $G_{\alpha }$ there are some $\beta _{1},\beta _{2} <\alpha $, with $G_{\alpha }(x)=G_{\beta _{1}}(x)$ and $G_{\alpha }(y)=G_{\beta _{2}}(y)$. The sequence is also decreasing, so for all $\beta _{1}\le \gamma \le \alpha $, $G_{\gamma }(x)=G_{\alpha }(x)$ and for all $\beta _{2}\le \gamma \le \alpha $, $G_{\gamma }(y)=G_{\alpha }(y)$. Let $\beta =max(\beta _{1},\beta _{2})$. Then $G_{\alpha }(x)=G_{\beta }(x)\neq G_{\beta }(y)=G_{\alpha }(y)$, since $G_{\beta }$ is a proper coloring. Thus, $G_{\alpha }$ is also a proper coloring and clearly $G_{\alpha }\le G_{\beta +1}=H(G_{\beta })<G_{\beta }$ for $\beta <\alpha $. The sequence $\langle G_{\alpha}:\alpha \in On \rangle$ is strictly decreasing so it is an injection from the class of ordinals into a set and that is a contradiction.
\end{proof}

\begin{lemma}\label{PRO_Thm5}
 A proper coloring $F:{}^{\lambda}\kappa \longrightarrow \mu $  is minimal 
 if and only if for all $x\in {}^{\lambda }{\kappa } $ and $\beta <F(x)$, there is $y\in {}^{\lambda }{\kappa } $ such that $y\perp x$ and $F(y)=\beta $.
\end{lemma}

\begin{proof}
 Assume that $F$ is minimal, $x\in {}^{\lambda }{\kappa } $ and $\beta <F(x)$. 
 Define $G$ as $G(y)=F(y)$ for $y\neq x$ and $G(x)=\beta $. Then clearly $G<F$ holds, so $G$ is not a proper coloring by the minimality of $F$.
 Since $G$ differs from $F$ only in one point, 
 there would be some $y\perp x$ with $F(y)=G(y)=G(x)=\beta$. So the condition holds.

Assume now that $F$ is not minimal, $G<F$ is a proper coloring.
Let 
\begin{displaymath}
{\beta}=\min\{G(x): G(x)<F(x)\}.
\end{displaymath}
Fix $x$ with $G(x)={\beta}$. Then ${\beta}=G(x)<F(x)$ and we claim that 
there is no $y\in {}^{\lambda }{\kappa } $ such that $y\perp x$ and $F(y)=\beta $.

Assume on the contrary that there is such a $y$.
Then $G(y)\le F(y)={\beta}$, so $G(y)={\beta}$ by the minimality 
${\beta}$. Thus, $G(y)=G(x)$ and $x\perp y$, meaning that $G$ is not a proper coloring. Contradiction, there is no suitable $y$, so the condition fails.
\end{proof}

Next we recall the definition of tightness from \cite{GLL}.

\begin{definition}\label{PRO_Def5}
 A proper coloring $F:{}^{\lambda}\kappa \longrightarrow \mu $  is \emph{tight} if for all $x\in {}^{\lambda }{\kappa } $ and $\beta \in \mu$ with $F(x)\neq \beta$, there is $y\in {}^{\lambda }{\kappa } $ with $y\perp x$ and $F(y)=\beta $.
 
 A proper coloring $F:{}^{\lambda}\kappa \longrightarrow \mu $  is \emph{$\nu $-tight} if for all sequences 
 $(x_{\alpha}:\alpha <\nu )\subset {}^{\lambda }{\kappa } $ and $\beta \in \mu$ with $F(x_{\alpha })\neq \beta$ for all $\alpha \in \nu$, there is $y\in {}^{\lambda }{\kappa } $ with $y\perp x_{\alpha }$ for all $\alpha \in \nu$ and $F(y)=\beta $. 
 
 A proper coloring $F:{}^{\lambda}\kappa \longrightarrow \mu $ is \emph{$<\nu $-tight} if it is $\eta $-tight for all $\eta <\nu $.
\end{definition}

Tightness of $F$ means that we cannot change $F$ in one point to remain proper. Clearly if $\eta <\nu $, then all $\nu $-tight proper colorings are $\eta $-tight. We have seen that minimality is a weak condition in the sense that  by Lemma  \ref{PRO_Lm4} we can find a minimal proper coloring under an arbitrary  proper coloring. By Lemma  \ref{PRO_Thm5} we have that tightness implies minimality. However, $<\kappa$-tightness  implies that $F$ is nice in some way.

\begin{lemma}\label{PRO_Lm5}
 If $F:{}^{\lambda}\kappa \longrightarrow \kappa $ is a ${<\kappa }$-tight proper coloring, then it is weakly uniform.
\end{lemma}

\begin{proof}
 We construct the functions $r_{\alpha }$ by induction. Choose any $r_{0}\in {}^{\lambda }{\kappa } $ with $F(r_{0})=0$. Suppose we have already chosen $r_{\beta }$-s for $\beta <\alpha $ that are pairwise totally different and $F(r_{\beta })=\beta $. Now use the $|\alpha |$-tightness condition for the $r_{\beta }$-s and the value $\alpha $. Clearly for all $\beta \in \alpha $, $F(r_{\beta })=\beta \neq \alpha $ holds, so we can choose $r_{\alpha }$, such that $r_{\alpha }\perp r_{\beta }$ for $\beta <\alpha $ and $F(r_{\alpha })=\alpha $. Then the conditions in Definition 3 hold for the $r_{\alpha }$-s.
\end{proof}

As we have seen, by Theorem \ref{PRO_Thm4}, weak uniformity gives strong restriction on $F$, even by Corollary \ref{PRO_Cor2}, if $\lambda $ is not too big (below first measurable cardinal), then $F$ depends only on one coordinate of the product. In the next section we will show that even 2-tightness gives a classification, and later we will see that tightness is a really weak property.

\section{The classification of 2-tight proper colorings}

In this section we will see that 2-tight proper colorings can be written in a "nice" form.

\begin{lemma}\label{PRO_Lm6}
 If $F:{}^{\lambda}\kappa \longrightarrow \mu $ is a 2-tight proper coloring, $x_{0},x_{1}\in {}^{\lambda }{\kappa } $ and $y\in {}^{\lambda }{\kappa } $ are such that $y(i)\in \{ x_{0}(i),x_{1}(i)\} $ for all $i\in \lambda $, then $F(y)\in \{ F(x_{0}),F(x_{1})\} $.
\end{lemma}

\begin{proof}
 Choose any $\beta \in \mu $ with $\beta \not\in \{ F(x_{0}),F(x_{1})\} $. Then by the 2-tightness of $F$, we can find a $z\in {}^{\lambda }{\kappa } $, such that $x_{0}\perp z\perp x_{1}$ and $F(z)=\beta $. Then for all $i\in \lambda $. $x_{0}(i)\neq z(i)\neq x_{1}(i)$, so $z(i)\neq y(i)$, thus $z\perp y$. Since $F$ is a proper coloring, we have $F(y)\neq F(z)=\beta $. Then the only possible values for $F(y)$ are $F(x_{0})$ and $F(x_{1})$.
\end{proof}

\begin{lemma}\label{PRO_Lm7}
 If $F:{}^{\lambda}\kappa \longrightarrow \mu $ is a 2-tight proper coloring, $x_{0},x_{1},...,x_{n-1}\in {}^{\lambda }{\kappa } $ for some $n\in \omega $ and $y\in {}^{\lambda }{\kappa } $ is such that $y(i)\in \{ x_{0}(i),...,x_{n-1}(i)\} $ for all $i\in \lambda$, then $F(y)\in \{ F(x_{0}),...,F(x_{n-1})\} $.
\end{lemma}

\begin{proof}
 We use induction on $n$. If $n=1$, the statement is trivial, if $n=2$, we can use Lemma  \ref{PRO_Lm6}. Suppose that it is true for some $n$ and let $x_{0},x_{1},...,x_{n}\in {}^{\lambda }{\kappa }$ and a $y\in {}^{\lambda }{\kappa } $ is such that $y(i)\in \{ x_{0}(i),...,x_{n}(i)\} $ for all $i\in \lambda$. Construct $z\in {}^{\lambda }{\kappa } $, such that $z(i)=y(i)$ if $y(i)\in \{ x_{0}(i),...,x_{n-1}(i)\} $ and $z(i)=x_{0}(i)$ otherwise. Then by the construction of $z$, we have $z(i)\in \{ x_{0}(i),...,x_{n-1}(i)\} $ for all $i\in \lambda $ so by the induction hypothesis, we have $F(z)\in \{ F(x_{0}),...,F(x_{n-1})\} $. On the other hand if $z(i)\neq y(i)$ for some $i \in \lambda$, then the only possible value for $y(i)$ is that $y(i)=x_{n}(i)$, so $y(i)\in \{ z(i),x_{n}(i) \}$ for all $i\in \lambda $. Then by Lemma  \ref{PRO_Lm6}, we have $F(y)\in \{ F(z),F(x_{n})\} \subseteq \{ F(x_{0}),...,F(x_{n})\} $.
\end{proof}

\begin{lemma}\label{PRO_Lm8}
  If $F:{}^{\lambda}\kappa \longrightarrow \mu $ is a 2-tight proper coloring, then there is an ultrafilter $\mathscr{U}$  on $\lambda $, such that for any sequence $x_{0},...,x_{n-1}\in {}^{\lambda }{\kappa } $ of  pairwise totally different  elements and for any $y\in {}^{\lambda }{\kappa } $ with $y(i)\in \{ x_{0}(i),...,x_{n-1}(i)\} $ for all $i\in \lambda $ and $k\in n$, we have $$F(y)=F(x_{k})\Leftrightarrow \{ i\in \lambda |y(i)=x_{k}(i) \} \in \mathscr{U}.$$
 \end{lemma}
 
 \begin{proof}
 
  If $a=\{x_0,\dots, x_{n-1}\}\subset {}^{\lambda}{\kappa}$ is a finite set,
  write 
  \begin{displaymath}
  R_a=\big\{y\in {}^{\lambda}{\kappa}\ | \ \forall i\in {\lambda}\ y(i)\in 
  \{x_0(i),.\dots, x_{n-1}(i)\}\ \big\}.
  \end{displaymath}

\begin{claim}\label{PRO_Lm:8-1}
  If 
  $a=\{x_0,\dots, x_{n-1}\}\subset {}^{\lambda}{\kappa}$ is a finite set of pairwise totally different functions, $n\ge 3$, then 
  there is an ultrafilter $\mathscr{U}_{a}$   on ${\lambda}$
 such that  
 \begin{displaymath}
 \tag{$\ddag_a$}
 \forall y\in R_a \ \big(
 F(y)=F(x_k)\ \Longleftrightarrow\ 
 \{i\in {\lambda}: y(i)=x_k(i)\}\in \mathscr{U}_a \ \big).
 \end{displaymath}
\end{claim}

\begin{proof}[Proof of the Claim]
There is a natural bijection $\varphi:{}^{\lambda}n\to R_a$ defined by the formula
$$\varphi(z)(i)=x_{z(i)}(i). $$
Since $F''R_a\subset \{F(x_0),\dots,F(x_{n-1})\}$ we can define 
$G:{}^{\lambda}n\to n$ such that 
\begin{displaymath}
\text{$G(z)=k$ iff $F(\varphi(z))=F(x_k)$.}
\end{displaymath}
Then $G$ is a proper coloring of ${}^{\lambda}n$ such that 
$G(c_k)=k$.

Then, by 
 Theorem \ref{PRO_TmKoTo}, using that $G(c_k)=k$ for $k<n$, there is an ultrafilter $\mathscr{U}_{a}$   on ${\lambda}$
 such that for each $z\in {}^{\lambda}n$, 
 \begin{displaymath}
 G(z)= k \ \Longleftrightarrow \ \{i\in {\lambda}: z(i)=k\}\in \mathscr{U}_a.
 \end{displaymath}
 Since $G(z)=k$ iff $F({\varphi}(z))=F(x_k)$
 and $\{i\in {\lambda}: z(i)=k\}=\{i\in {\lambda}: \varphi(z)(i)=x_k(i)\}$
, we verified $(\ddag_a)$.
\end{proof}

\begin{claim}\label{PRO_Lm:8-2}
  If 
  $a\subset b$ are finite sets of pairwise totally different functions, $|a|\ge 3$, then $\mathscr{U}_a=\mathscr{U}_b$. 
\end{claim}

\begin{proof}[Proof of the Claim]
Let $A\in \mathscr{U}_a$ be arbitrary.
  Pick $\{x_0,x_1\}\in {[a]}^{2}$, and 
  fix $x\in {}^{\lambda}{\kappa}$ such that 
  $x\equiv_A x_0$ and $x\equiv_{{\lambda}\setminus A}x_1$.
  Then $F(x)=F(x_0)$ by ($\ddag_a$).
  But ($\ddag_{b}$) holds  and so  $F(x)=F(x_0)$ implies that 
  $A=\{i: x(i)=x_k(i)\}\in \mathscr{U}_{b}$.
  Thus, $\mathscr{U}_{a}\subset \mathscr{U}_{b}$, 
  but both are ultrafilters on $\lambda $, so 
  $\mathscr{U}_{a}=\mathscr{U}_{b}$.
\end{proof}

\begin{claim}\label{PRO_Lm:8-3}
  If 
  $a$ and $b$ are finite sets of pairwise totally different functions, 
  $|a|,|b|\ge 3$, then $\mathscr{U}_a=\mathscr{U}_b$. 
  (We do not assume that $x\perp y$ for $x\in a$ and $y\in b!$)
\end{claim}
\begin{proof}
Pick a 3-element set $c=\{c_0,c_1,c_2\}\in {[{}^{\lambda}{\kappa}]}^{3}$
of pairwise totally different functions 
such that $c_\ell(i)\notin\{x(i):x\in a\cup b\}$. 
Then the elements of $a\cup c$ and $b\cup c$ are pairwise totally 
different, so $\mathscr{U}_{a\cup c}$ and  $\mathscr{U}_{b\cup c}$
 are defined, and, by Claim \ref{PRO_Lm:8-2}, 
 \begin{displaymath}
 \mathscr{U}_{a}=\mathscr{U}_{a\cup c}=\mathscr{U}_{c}=
 \mathscr{U}_{b\cup c}=\mathscr{U}_{b}
 \end{displaymath}
 which was to be proved. 
\end{proof}

We are ready to complete the proof of Lemma \ref{PRO_Lm8}.
By Claim \ref{PRO_Lm:8-3} there is an ultrafilter $\mathscr{U}$ on ${\lambda}$
such that 
$\mathscr{U}_a=\mathscr{U}$
whenever  
$a$ is a finite set of pairwise totally different functions
with  $|a|\ge 3$.

Let $m=n+2$, and pick $x_n, x_{n+1}$ such that 
$a=\{x_0,\dots, x_{m-1}\}\subset {}^{\lambda}{\kappa}$ is a finite set of pairwise totally different functions. Then $|a|=n+2\ge 3$, so
we have  $\mathscr{U}_a=\mathscr{U}$,
i.e. 
\begin{displaymath}
  \tag{$\ddag$}
  \forall y\in R_a \ \big(
  F(y)=F(x_k)\ \Longleftrightarrow\ 
  \{i\in {\lambda}: y(i)=z_k(i)\}\in \mathscr{U} \ \big).
  \end{displaymath}
 \end{proof}

\begin{definition}\label{PRO_Def6}
 For a 2-tight proper coloring $F:{}^{\lambda}\kappa \longrightarrow \mu $  we define the \emph{corresponding ultrafilter} as the unique  ultrafilter $\mathscr{U}$ on $\lambda $, such that for every finite set  $x_{0},...,x_{n-1}\in {}^{\lambda }{\kappa } $ of pairwise totally different function and $y\in {}^{\lambda }{\kappa } $ with $y(i)\in \{ x_{0}(i),...,x_{n-1}(i)\} $ for all $i\in \lambda $ and $k\in n$, we have 
 $$F(y)=F(x_{k})\ \Leftrightarrow\ \{ i\in \lambda |y(i)=x_{k}(i) \} \in \mathscr{U}.$$ 
 By Lemma  \ref{PRO_Lm8} such an ultrafilter exists and clearly unique.
\end{definition}

\begin{lemma}\label{PRO_Lm9}
 If $F:{}^{\lambda}\kappa \longrightarrow \mu $ is a 2-tight proper coloring, $\mathscr{U}$ is the corresponding ultrafilter and $x,y\in {}^{\lambda }{\kappa } $, then $$F(x)=F(y)\ \Longleftrightarrow\ \{ i\in \lambda |x(i)=y(i)\} \in \mathscr{U}.$$
\end{lemma}

\begin{proof}
 Let $z\in {}^{\lambda }{\kappa } $ be such that $z(i)=y(i)$ if $x(i)\neq y(i)$ and $z(i)\neq x(i)$ if $x(i)=y(i)$. Then $z\perp x$, so $F(z)\neq F(x)$. Also we defined $z$ as $y(i)\in \{ x(i),z(i)\}$ holds for all $i\in \lambda $, so we can use Lemma  \ref{PRO_Lm8}, for $x,z$ and $y$, so we have $F(x)=F(y)\Leftrightarrow \{ i\in \lambda |x(i)=y(i)\} \in \mathscr{U}$, and that is what we wanted to prove.
\end{proof}

\begin{definition}\label{PRO_Def7}
 For a fixed ultrafilter $\mathscr{U}$  on $\lambda $ and $x\in {}^{\lambda }{\kappa } $ we define its equivalence class 
 $\eclass{x}{U}=\{ x'\in {}^{\lambda }{\kappa } |\{i\in \lambda |x'(i)=x(i)\} \in \mathscr{U} \}$.
 
 We define the quotient $\bigslant{{}^{\lambda }{\kappa } }{\mathscr{U}}=\{ \eclass{x}{U}: x\in {}^{\lambda }{\kappa } \}$. For an $X\in \bigslant{{}^{\lambda }{\kappa } }{\mathscr{U}}$, we say an $r\in {}^{\lambda }{\kappa } $ is a representative for $X$ if $r\in X$ or equivalently $\eclass{r}{U}=X$.
\end{definition}

\begin{theorem}\label{PRO_Thm6}
Given a function   $F:{}^{\lambda}\kappa \longrightarrow \mu $
the following statements are equivalent
\begin{enumerate}[(1)]
\item $F$ is a 2-tight proper coloring, 
\item $F$ is a $k$-tight proper coloring   for each $k<{\omega}$,
\item there is an ultrafilter $\mathscr{U}$ on $\lambda $ and a bijection $h:\bigslant{{}^{\lambda }{\kappa } }{\mathscr{U}} \longrightarrow \mu$, such that $F(x)=h(\uclass{x})$ for all $x\in {}^{\lambda }{\kappa } $.
\end{enumerate}
\end {theorem}

\begin{proof}(1) $\longrightarrow$ (3):
 Since $F$ is 2-tight, we can consider 
   the corresponding ultrafilter $\mathscr{U}$. 
 By Lemma \ref{PRO_Lm9}, the following definition is meaningful:
\begin{displaymath}
h(\uclass{x})={\beta}\ \text{ iff } \ F(x')={\beta}
\text{ for each }x'\in \uclass{x}.
\end{displaymath}
 
If $\uclass{x}\ne \uclass{y}$, then $F(x)\ne F(y)$ by 
Lemma \ref{PRO_Lm9}, so $h(\uclass{x})=F(x)\ne F(y)=h(\uclass{y})$. Thus, 
$h$ is injective.

Since $F$ is onto, so is $h$. Thus, $h$ has the required properties.



\medskip
\noindent (3) $\longrightarrow$ (2)

First we will show that $F$ is a proper coloring. If $x,y\in {}^{\lambda }{\kappa } $ with $x\perp y$, then $\{ i\in \lambda |x(i)=y(i)\} =\emptyset \not\in \mathscr{U}$, thus $\uclass{x}\neq \uclass{y}$. Since $h$ is a bijection, we have $F(x)=h(\uclass{x})\neq h(\uclass{y})=F(y)$. 

Now we need to show that $F$ is $k$-tight. Let $x_{0},\dots,x_{k-1}\in {}^{\lambda }{\kappa } $ be arbitrary and $\beta \in \mu $ with $\beta \not\in 
\{ F(x_{0}),\dots, F(x_{k-1})\} $. Since $h$ is a bijection, there is a 
$z\in {}^{\lambda}{\kappa}$ with $h(\uclass{z})={\beta}$.
Then $F(z)={\beta}$ as well. 

Let $A_{\ell}=\{ i\in \lambda : z(i)\neq x_{\ell}(i)\} $ for $\ell<k$.
Since $h(\uclass{x_\ell})=F(x_\ell)\ne {\beta}=h(\uclass{z})$,
it follows that $A_\ell\in \mathscr{U}$ by Lemma \ref{PRO_Lm9}.
Thus, $A=\bigcap_{\ell<k}A_\ell\in \mathscr{U}$ as well. 


Let $y\in {}^{\lambda }{\kappa } $ be such that $y\equiv _{A}z$ and 
$y(i)\in {\kappa}\setminus \{x_0(i),\dots, x_{k-1}(i)\}$
for $i\in {\lambda}\setminus A$.

Then $x_\ell\perp y$ for each $\ell<k$, but 
$\{i\in {\lambda}: y(i)=z(i)\}\in \mathscr{U}$, and so 
$F(y)=h(\uclass{y})=h(\uclass{z})=F(z)={\beta}$.

So $F$ is $k$-tight. 

\medskip
\noindent (2) $\longrightarrow$ (1). Special case. 
\end{proof}

\begin{corollary}\label{PRO_cr:t2}
If ${\kappa},{\mu}<2^{\omega}$ are infinite cardinals, and $
F:{}^{\omega}{\kappa}\to {\mu}$ is a 2-tight proper coloring, then 
${\kappa}={\mu}$ and there are $i<{\omega}$
and a permutation ${\pi}\in Sym({\kappa})$ such that 
$F(x)={\pi}(x(i))$ for each $x\in {}^{\lambda}{\kappa}$.
\end{corollary}

\begin{proof}
By theorem \ref{PRO_Thm6} there is an ultrafilter $\mathscr{U}$ on ${\omega}$ and a bijection $h:\bigslant{{}^{\omega }{\kappa } }{\mathscr{U}} \longrightarrow \mu$, such that $F(x)=h(\uclass{x})$ for all $x\in {}^{\omega}\kappa $.

  If $\mathscr U$ is a non-principal ultrafilter on ${\omega}$, then 
$|{}^{\omega}{\kappa}/\mathscr U|=2^{\omega}>\mu $.  Thus, $\mathscr U$ is 
a principal ultrafilter: there is $i\in {\omega}$ such that 
$\mathscr U=\{U\subset {\omega}:i\in U\}$. 
Then $F(x)=h({\{y\in {}^{\omega}{\kappa}:y(i)=x(i)\}})$. 
So if we define ${\pi}\in Sym({\kappa})$ by the formula
\begin{displaymath}
{\pi}({\alpha})=h(\{y\in {}^{\omega}{\kappa}:y(i)={\alpha}\}),
\end{displaymath}
then $F(x)={\pi}(x(i))$ for each $x\in {}^{\omega}{\kappa}$.
\end{proof}

This way we can characterize all 2-tight proper colorings: they can be  constructed from an ultrafilter $\mathscr{U}$ on $\lambda $ and a bijection $h:\bigslant{{}^{\lambda }{\kappa } }{\mathscr{U}} \longrightarrow \mu $. To have such a bijection, we need that $|\bigslant{{}^{\lambda }{\kappa } }{\mathscr{U}}|=\mu $. For example, if $\mathscr{U}$ is principal, then the quotient has cardinality $\kappa $ (this case gives the trivial colorings). However, in general it is hard to compute $|\bigslant{{}^{\lambda }{\kappa } }{\mathscr{U}}|$ for an arbitrary ultrafilter.




\begin{corollary}\label{PRO_cr:t3}
If $F:{}^{\omega }{\omega }\longrightarrow \omega $ is a 2-tight proper coloring, then $F$ is trivial.
\end{corollary}


As a final theorem of this section, we characterize the $\nu $-tight proper colorings for an infinite cardinal $\nu $.

\begin{theorem}\label{PRO_Thm7}
 A function  $F:{}^{\lambda}\kappa \longrightarrow \mu $ is a $\nu $-tight proper coloring for an infinite cardinal $\nu <\kappa $  if and only if there is a $\nu ^{+}$-complete ultrafilter $\mathscr{U}$ on $\lambda $ and a bijection $h:\bigslant{{}^{\lambda }{\kappa } }{\mathscr{U}} \longrightarrow \mu$, such that $F(x)=h(\uclass{x})$ for all $x\in {}^{\lambda }{\kappa } $.
\end{theorem}

\begin{proof}
 First suppose that $F$ is a $\nu $-tight proper coloring. Then it is 2-tight, so by Theorem \ref{PRO_Thm6}, there is an ultrafilter $\mathscr{U}$ on $\lambda $ and a bijection $h:\bigslant{{}^{\lambda }{\kappa } }{\mathscr{U}} \longrightarrow \mu$, such that $F(x)=h(\uclass{x})$ for all $x\in {}^{\lambda }{\kappa } $. We only need to check that $\mathscr{U}$ is $\nu ^{+}$-complete. Suppose for contradiction that $\mathscr{U}$ is not $\nu ^{+} $-complete and pick $A_{\alpha}\in \mathscr{U}$-s for $\alpha \in \nu $ with $A=\bigcap_{\alpha \in \nu }\not\in \mathscr{U}$. Then for each ${\alpha}<{\nu}$ consider  
 the  element  $x_{\alpha }$ of $^{\lambda }{\kappa }$, such that $x_{\alpha }\equiv _{A_{\alpha }}c_{1}$ and $x_{\alpha }\equiv _{\lambda -A_{\alpha }}c_{0}$ for all $\alpha \in \nu $. Then $\{i\in \lambda |x_{\alpha }(i)=c_{1}(i)\} =A_{\alpha }\in \mathscr{U}$, so $\uclass{x_{\alpha }}=\uclass{c_{1}}$, thus $F(x_{\alpha })=h(\uclass{x_{\alpha }})=h(\uclass{c_{1}})=F(c_{1})$ for all $\alpha \in \nu $. Clearly $c_{0}\perp c_{1}$, so $F(c_{0})\neq F(c_{1})$. This means that $F(c_{0})\not\in \{ F(c_{1})\} =\{ F(x_{\alpha }):\alpha \in \nu \}$, so by $\nu $-tightness, there is a $y\in {}^{\lambda }{\kappa } $ with $y\perp x_{\alpha }$ for $\alpha \in \nu $ and $F(y)=F(c_{0})$. Then $h(\uclass{y})=h(\uclass{c_{0}})$, so $\uclass{y}=\uclass{c_{0}}$, thus $B= \{i\in \lambda |y(i)=0\} \in \mathscr{U}$. For any $i\in B$ and $\alpha \in \nu $, since $y\perp x_{\alpha }$ and $y(i)=0$, we have $x_{\alpha }(i)\neq 0$, but then $i\in A_{\alpha }$. It is true for all $\alpha \in \nu $, so $i\in \bigcap_{\alpha \in \nu }A_{\alpha }=A$, thus $B\subseteq A$. This is a contradiction, since $A\not\in \mathscr{U},B\in \mathscr{U}$.
 
 Now suppose that $F$ is in that form, and we need to show that this is a $\nu$-tight proper coloring. By Theorem \ref{PRO_Thm6},  $F$ is a proper coloring, so we need to show that $F$ is $\nu $-tight. Let $\{ x_{\alpha }:\alpha \in \nu \}$ be elements of ${}^{\lambda }{\kappa } $ and $\beta \in \mu $ with $\beta \not\in \{F(x_{\alpha }):\alpha \in \nu \}$. Let 
 $z\in {}^{\lambda }{\kappa }$ with $h(\uclass{z})=\beta $. 
 For all $\alpha \in \nu $ let $A_{\alpha }=\{i\in \lambda |x_{\alpha }(i)\neq z(i)\}$. 
 Since $h(\uclass{x_{\alpha}})=F(x_{\alpha})\ne {\beta}=h(\uclass{z})$,
it follows that $A_{\alpha}\in \mathscr{U}$ by Lemma \ref{PRO_Lm9}.
 Since $\mathscr{U}$ is $\nu ^{+}$-complete, $A=\bigcap _{\alpha \in \nu }A_{\alpha}\in \mathscr{U}$. Construct a $y\in {}^{\lambda }{\kappa } $ with $y\equiv _{A}z$ and $y\perp _{\lambda -A}x_{\alpha }$ for every $\alpha \in \nu $. We can construct such $y$, since for any $i\in \lambda -A$ $|\{x_{\alpha }(i):\alpha \in \lambda \} |\le \nu <\kappa $, so we can choose $y(i)$ that differs from all. Then $y\perp x_{\alpha}$ for all $\alpha \in \nu $, and $\{i\in \lambda |y(i)=z(i)\} =A\in \mathscr{U}$, thus $\uclass{y}=\uclass{{z}}$, so $F(y)=h(\uclass{y})=h(\uclass{z})=\beta $. Then $y$ satisfies the requirements, so $F$ is $\nu $-tight.
\end{proof}

\begin{corollary}\label{PRO_Cor4}
 If $\lambda $ is below the smallest measurable cardinal and $F:{}^{\lambda}\kappa \longrightarrow \mu $ is an $\omega $-tight proper coloring then $\mu =\kappa $ and $F(x)=\pi (x(i))$ for some $i\in \lambda$ and $\pi \in Sym(\kappa )$.
\end{corollary}

\begin{proof}
 By Theorem \ref{PRO_Thm7}, there  is an $\omega _{1}$-complete ultrafilter $\mathscr{U}$ on $\lambda $ and a bijection $h:\bigslant{{}^{\lambda }{\kappa } }{\mathscr{U}} \longrightarrow \mu$, such that $F(x)=h(\uclass{x})$ for all $x\in {}^{\lambda }{\kappa } $. Since $\lambda $ is below the first measurable cardinal, $\mathscr{U}$ must be principal, so there is an $i\in \lambda $, such that $\mathscr{U}=\{A\subseteq \lambda |i\in A\}$. For all $x\in {}^{\lambda }{\kappa } $, clearly $x(i)=c_{x(i)}(i)$, so $\{j\in \lambda |x(j)=c_{x(i)}(j)\} \in \mathscr{U}$.
 On one hand it means that $\bigslant{{}^{\lambda }{\kappa } }{\mathscr{U}}=\{ \uclass{c_{\alpha }}:\alpha \in \kappa \}$, so $\mu =|\bigslant{{}^{\lambda }{\kappa } }{\mathscr{U}}|=\kappa $. We can define $\pi (\alpha )=h(\uclass{c_{\alpha }})=F(c_{\alpha })$. Since $h$ is a bijection, $\pi $ is also, so $\pi \in Sym(\kappa )$. On the other hand for arbitrary $x$ we have $F(x)=h(\uclass{x})=h(\{\uclass{c_{x(i)}})=\pi (x(i))$. 
\end{proof}

\section{Finitely independent proper colorings}

In this section, we try to give a bunch of proper colorings that are far away from trivial.

\begin{definition}\label{PRO_Def8}
For $x,y\in {}^{\lambda }{\kappa } $ we say $x$ and $y$ are\emph{ almost equal}, denoted by $x{\equiv }^{*}y$, if $|\Delta(x,y)|<\omega $, and \emph{almost totally different}, denoted by $x\perp ^{*}y$ if $|\lambda -\Delta(x,y)|<\omega $. 

A proper coloring $F:{}^{\lambda}\kappa \longrightarrow \mu $  is \emph{finitely independent} if $F(x)=F(y)$ whenever $x{\equiv }^{*}y$. $F$ is \emph{strongly proper} if $F(x)\neq F(y)$ whenever $x\perp ^{*}y$.
\end{definition}

Clearly strongly proper colorings are proper, and finitely independent proper colorings are not trivial. 

\begin{lemma}\label{PRO_Lm10}
A  finitely independent proper coloring $F:{}^{\lambda}\kappa \longrightarrow \mu $ is   strongly proper.
\end{lemma}

\begin{proof}
 Let $x,y\in {}^{\lambda }{\kappa } $ with $x\perp ^{*}y$. Choose $z\in {}^{\lambda }{\kappa }$ with $z\equiv _{\Delta (x,y)}y$ and $z\perp _{\lambda -\Delta (x,y)}x$. Then $|\Delta (y,z)|=|\lambda -\Delta (x,y)|<\omega $, so $y{\equiv }^{*}z$, and since $F$ is finitely independent $F(z)=F(y)$. On the other hand $z\perp x$, so $F(x)\neq F(z)=F(y)$.
\end{proof}

We have already seen in Corollary \ref{PRO_cr:t2} that all 2-tight proper colorings ${F:{}^{\omega }\omega \longrightarrow \omega }$ are trivial, and if CH does not hold, then all 2-tight colorings ${F:{}^{\omega }\omega _{1} \longrightarrow \omega _{1}}$ are trivial, as well, thus they cannot be finitely independent. On the other hand, by  Theorem \ref{PRO_Thm7} we can see that all non-trivial 2-tight proper colorings must be finitely independent. It is natural to ask if there are finitely independent, tight proper colorings $F:{}^{\omega }\omega \longrightarrow \omega $ or $F:{}^{\omega }\omega _{1} \longrightarrow \omega _{1}$

The statement of the following lemma is well-known.
\begin{lemma}\label{PRO_Lm11}
 There is a set $\mathscr{A}\subset {}^{\omega }\omega $ with $|\mathscr{A}|=2^{\omega }$ such that $x\perp ^{*}y$ for all $x,y\in \mathscr{A}$ with $x\neq y$.
\end{lemma}


\begin{corollary}\label{PRO_Cor5}

If ${\omega}\le {\kappa}\le 2^{\omega}$, then 
\begin{displaymath}
2^{\omega}=\min\{{\mu}:\text{there is a finitely independent proper coloring 
$F:{}^{\omega}{\kappa}\to {\mu}\}$}
\end{displaymath}
Hence there is no finitely independent proper coloring $F:{}^{\omega }\omega \longrightarrow \omega $ and if CH does not hold, then there is no finitely independent proper coloring $F:{}^{\omega }\omega _{1}\longrightarrow \omega _{1}$.

\end{corollary}

\begin{proof}
 Suppose $F:{}^{\omega}{\kappa}\to {\mu}$ is a finitely independent proper coloring.  Then, by Lemma  \ref{PRO_Lm10}, $F$ is strongly proper. Let $\mathscr{A}$ be the set from Lemma  \ref{PRO_Lm11}.
 Then $F|_\mathscr{A}$ must be one-to-one, so $\mu \ge |Ran(F)|\ge |\mathscr{A}|=2^{\omega }$.

 Now let $\mathscr{U}$ be a non-principal ultrafilter on $\omega $, then  $|\bigslant{^{\omega }\kappa }{\mathscr{U}}|\le |^{\omega }\kappa |=2^{\omega }$. 
 Let $h: \bigslant{^{\omega }\kappa }{\mathscr{U}} \longrightarrow 2^{\omega }$ be a bijection, so by theorem \ref{PRO_Thm6}, the function $F:{}^{\omega }\kappa \longrightarrow 2^{\omega }$ defined by $F(x)=h(x_{\mathscr{U}})$ is a proper coloring, which is even 2-tight, and it is finitely independent, because for any $x,y\in {}^{\omega }\kappa $ with $|\Delta (x,y)|<\omega $ we have $x_{\mathscr{U}}=y_{\mathscr{U}}$, since $\mathscr{U}$ is non-principal, so $F(x)=h(x_{\mathscr{U}})=h(y_{\mathscr{U}})=F(y)$.
\end{proof}

Now we have seen that the existence of finitely independent tight proper colorings on $^{\omega }\omega _{1}$ to $\omega _{1}$ depends on CH. Now we will show that using CH we can construct a wide variety of tight finitely independent proper colorings. For this we will use notions of forcing. 

\begin{definition}\label{PRO_Def9}
 We define the {\em tight-coloring poset} $P$. The elements $p\in P$ are functions with $Dom(p)\subset {}^{\omega }\omega _{1}$, $Ran(p)\subset \omega _{1}$ such that  the following conditions hold:
 
\begin{enumerate}[1.]
\item 
 $|Dom(p)|\le \omega$,
\item 
there are no $x\neq y\in Dom(p)$ with $x{\equiv }^{*}y$,
\item  
  if $x,y\in Dom(p)$ and $x\perp ^{*}y$, then $p(x)\neq p(y)$.
\end{enumerate}
(Observe that  $\neg (x\equiv^* y)$ does not imply $x\perp^* y$ .)

 We define $q\le p$ if $p\subseteq q$
\end{definition}

Clearly $\langle P,\le \rangle$ is a poset and the empty function is the greatest element of $P$.

\begin{lemma}\label{PRO_Lm12}
 The tight-coloring $P$ is $\omega $-closed.
\end{lemma}

\begin{proof}

Straightforward.
\end{proof}

\begin{lemma}\label{PRO_Lm13}
 If $x\in {}^{\omega }\omega _{1}$ then 
 $$D_{x}=\{ p\in P\ |\ \exists x'\in Dom(p)\ x'{\equiv }^{*}x\}$$ is an open dense subset of $P$.
\end{lemma}

\begin{proof}Straightforward.
\end{proof}

\begin{lemma}\label{PRO_Lm14}
 If $x\in {}^{\omega }\omega _{1}$ and $\alpha \in \omega _{1}$, then 
 \begin{displaymath}
  E_{x,\alpha }=\{ p\in P\ |\ \exists x'\in Dom(p)[(x'\equiv^* x\lor x'\perp^* x)\land p(x')={\alpha}]\}.
 \end{displaymath}
 is an open dense subset of $P$.
\end{lemma}

\begin{proof}
 Clearly $E_{x,\alpha }$ is open. Let $p\in P$ be arbitrary. By Lemma  \ref{PRO_Lm13} $D_{x}$ is open dense, so we can choose $q\in D_{x}$ with $q\le p$. Let $x'\in {}^{\omega }\omega _{1}$ be the element such that $x'{\equiv }^{*}x$ and $x'\in Dom(q)$. If $q(x')=\alpha $, then we are done, so suppose  $q(x')\neq \alpha $. Let $(t_{n})_{n\in \omega}$ be an enumeration of $Dom(q)$ and $z_{n}(n\in \omega ) $ be an enumeration of $\{ y\in Dom(q)|q(y)=\alpha \}$ (if any of these sets are finite we can repeat the elements, but all elements must be listed). Also let $S:\omega \longrightarrow \omega \times 3$ be a function that all elements of $\omega \times 3$ have infinitely many preimages. We will define a sequence $(f_{k})_{k\in \omega }$ of functions, such $f_{0}\subseteq f_{1}\subseteq f_{2}\subseteq ...$, $Dom(f_{k})\subset \omega$ is finite and $Ran(f_{k})\subset \omega _{1}$.
 
 The sequence is defined by induction, let $f_{0}=\emptyset $. If $S(k)=(n,0)$ for some $n\in \omega $, and $n\not\in Dom(f_{k})$, then choose any $\beta \neq x'(n)$ and let $f_{k+1}=f_{k}\cup \{ (n,\beta )\} $, if $n\in Dom(f_{k})$, then simply $f_{k+1}=f_{k}$. In the case $S(k)=(n,1)$ choose any $m_{k}\not\in Dom(f_{k})$  and $\beta \in \omega _{1}$ with $x'(m_{k})\neq \beta \neq t_{n}(m_{k})$ and let $f_{k+1}=f_{k}\cup \{ (m_{k},\beta )\}$. In the case $S(k)=(n,2)$ choose $m_{k}\not\in Dom(f_{k})$ such that $z_{n}(m_{k})\neq x'(m_{k})$. We can choose such $m_{k}$ since by $f(x')\neq \alpha =f(z_{n})$ we have $x~\neq z_{n}$ and they are in $Dom(q)$, so they differ in infinitely many places, and $Dom(f_{k})$ is finite. Let $f_{k+1}=f_{k}\cup \{ (m_{k},z_{n}(m_{k}))\} $.
 
 Now let $y=\bigcup _{k\in \omega }f_{k}$. It's the union of a chain of functions, so $y$ is clearly a function with $Dom(y)\subseteq \omega $, and $Ran(y)\subseteq \omega_{1}$. We need to show that $Dom(y)=\omega $. Let $n\in \omega $, then we can choose $k$ such that $S(k)=(n,0)$. Then by definition of $f_{k+1}$ we have $n\in Dom(f_{k+1})\subseteq Dom(y)$. Thus, $y\in {}^{\omega }\omega _{1}$. Let $q'$ be the function with $Dom(q')=Dom(q)\cup \{ y \}$, $q\subseteq q'$ and $q'(y)=\alpha $. We need to show that $q'\in P$. Property 1 is clear, as we defined it in only 1 extra point. Property 2 would be violated only if there was a $t\in Dom(q)$ with $t{\equiv }^{*}y$, but then $t=t_{n}$ for some $n\in \omega $, and there are infinitely many $k$-s with $S(k)=(n,1)$. For all of those $k$-s $f_{k+1}$ is defined as $y(m_{k})=f_{k+1}(m_{k})\neq t_{n}(m_{k})$, so they differ in infinitely many places. Property 3 would be violated by some $z\in Dom(q)$ if $z\perp ^{*}y$ and $q'(z)=q(z)=\alpha $. But then $z=z_{n}$ for some $n\in \omega $, and there are infinitely many $k$-s with $S(k)=(n,2)$. For all of those $k$-s we defined $f_{k+1}$ as $y(m_{k})=f_{k+1}(m_{k})=z_{n}(m_{k})$, so $y$ and $z_{n}$ agrees on infinitely many places, and that is a contradiction. So $q'\in P$ holds, clearly $q'\le q\le p$. Throughout the construction, $f_{k}$-s, we always maintain that for all $m\in Dom(f_{k})$ we have $f_{k}(m)\neq x'(m)$, thus $y\perp x'$. Then $x',y\in Dom(q')$, $x'{\equiv }^{*}x$, $y\perp ^{*} x$ and $q'(y)=\alpha $, so $q'\in E_{x,\alpha }$, thus $E_{x,\alpha }$ is dense. 
\end{proof}

Now, to continue, we need the following well-known theorem:

\begin{theorem}\label{PRO_Lm15} 
(Sikorski)
 If $P$ is an $\omega $-closed notion of forcing, $p\in P$, and $D_{\alpha }\subseteq P$ is a dense subset for $\alpha \in \omega _{1}$ , then there is a filter $G\subseteq P$ such that $p\in G$ and $G\cap D_{\alpha }\neq \emptyset $ for all $\alpha <\omega_{1}$. 
 
 
 
 
 
\end{theorem}

\begin{lemma}\label{PRO_Lm16}
 If $G$ is any filter on $P$ then $F_{0}=\bigcup G$ is a function.
\end{lemma}

\begin{proof}
Straightforward.
\end{proof}

\begin{lemma}\label{PRO_Lm17}

 If $G$ is a filter on $P$, such that $G\cap D_{x}\neq \emptyset $ for each $x\in {}^{\omega }\omega _{1}$  and we take $F_{0}=\bigcup G$, then for all $x\in {}^{\omega }\omega _{1}$, there is a unique $x'\in Dom(F_{0})$ with $x{\equiv }^{*}x'$.
\end{lemma}

\begin{proof}
Straightforward.
\end{proof}

\begin{lemma}\label{PRO_Lm18}
 If $G$ is a filter on $P$ such that for all $x\in {}^{\omega}\omega _{1}$ and $\alpha \in \omega _{1}$ we have $G\cap D_{x}\neq \emptyset $ and $G\cap E_{x,\alpha }\neq \emptyset $, and we take $F_{0}=\bigcup G$, then there is an  $F: {}^{\omega } \omega _{1}\longrightarrow \omega_{1}$ finitely independent, tight proper coloring, such that $F_{0}\subseteq F$. 
\end{lemma}

\begin{proof}
 By Lemma \ref{PRO_Lm17}, for all $x\in {}^{\omega }\omega _{1}$, there is a unique $x'\in Dom(F_{0})$ that $x{\equiv }^{*}x'$. Using this, we can define $F$ as $F(x)=F_{0}(x')$, and this is a function from $^{\omega }\omega _{1}$. Clearly if $x\in Dom(F_{0})$, then $x'=x$, so $F(x)=F_{0}(x)$, $F_{0}\subseteq F$ holds.
 
 Now we need to prove $F$ is a finitely independent, tight, proper coloring. First, we show that it is a proper coloring. Let $x,y\in {}^{\omega }\omega _{1}$ with $x\perp y$. Then there are $x',y'\in Dom(F_{0})$ with $x{\equiv }^{*}x'$, $y{\equiv }^{*}y'$ and $F(x)=F_{0}(x')$, $F(y)=F_{0}(y')$. Then $x'\perp ^{*}y'$ holds. Let $p,q\in G$ such that $x'\in Dom(p)$, $y'\in Dom(q)$ and choose $r\in G$ with $r\le p,q$. Then $y',x'\in Dom(r)$ and by Property 3 on $r$, we have $F(x)=F_{0}(x')=r(x')\neq r(y')=F_{0}(y')=F(y)$, so $F$ is a proper coloring.
 
 For the finite independence, choose $x,y\in {}^{\omega }\omega _{1}$ with $x{\equiv }^{*}y$. Let $x'\in Dom(F_{0})$ be such that $x{\equiv }^{*}x'$. But then also $y{\equiv }^{*}x'$, so $F(x)=F_{0}(x')=F(y)$.
 
 The last thing we need to prove is that $F$ is tight. Choose any $x\in {}^{\omega }\omega _{1}$ and $\alpha \in \omega _{1}$ with $F(x)\neq \alpha $ and $p\in G\cap E_{x,\alpha }$. Then there is an $x'\in Dom(p)$ with $x{\equiv }^{*}x'$. Clearly $x'\in Dom(F_{0})$, so $F(x)=F_{0}(x')$ is defined by this $x'$. Moreover either $p(x')=\alpha $ or there is a $y\in Dom(p)$ with $y\perp ^{*}x'$ and $p(y)=\alpha $. In the former case $F(x)=F_{0}(x')=p(x')=\alpha $, witch contradicts our assumption, so the latter case holds. Then $x\perp ^{*}y$ also holds and $F_{0}(y)=p(y)=\alpha $. Let $z\in {}^{\omega }\omega _{1}$ be such that $z(i)=y(i)$ if $y(i)\neq x(i)$ and $z(i)\neq x(i)$ if $y(i)=x(i)$. Then $x\perp z$ holds and $\Delta (z,y)=\omega -\Delta (x,y)$, finite, so $y{\equiv }^{*}z$ means $F(z)=F_{0}(y)=\alpha $, thus $F$ is tight. 
\end{proof}

Now we formulate our results.

\begin{theorem}\label{PRO_Thm8}
 If CH holds, $\mathscr{A}\subseteq {}^{\omega }\omega _{1}$, $|\mathscr{A}|\le \omega $ and  $p:\mathscr{A}\longrightarrow \omega _{1}$ is a partially finitely independent and partially strongly proper function (i.e. for all $x,y\in \mathscr{A}$ if $x{\equiv }^{*}y$, then $p(x)=p(y)$, and if $x\perp ^{*}y$ then $p(x)\neq p(y)$), then there is a finitely independent, tight proper coloring $F:{}^{\omega }\omega _{1}\longrightarrow \omega _{1}$ with $p\subset F$.
\end{theorem}

\begin{proof}
 Choose $\mathscr{B}\subseteq \mathscr{A}$, such that for all $x\in \mathscr{A}$, there is a unique $x'\in \mathscr{B}$ with $x\equiv ^{*}x'$. Let $q=p|_{\mathscr{B}}$. We show that $q\in P$, where $P$ is the tight-coloring. Property 1 clearly holds, as $|Dom(q)|=|\mathscr{B}|\le |\mathscr{A}|\le \omega $. For property 2 if $x,y\in Dom(q)=\mathscr{B}$, $x\neq y$, then $x\equiv ^{*}{y}$ is not possible by the choice of $\mathscr{B}$. For property 3 if $x,y\in Dom(q)$ with $x\perp ^{*}y$, then since $p$ is partially strongly proper, we have $q(x)=p(x)\neq q(y)=p(y)$. By CH, we have $|\{D_{x}:x\in {}^{\omega }\omega _{1} \} \cup \{E_{x,\alpha }:x\in {}^{\omega }\omega _{1} , \alpha \in \omega _{1}\}|=\omega _{1}$, so by Lemma  \ref{PRO_Lm15} there is a filter $G\subseteq P$, such that for all $x\in {}^{\omega }\omega _{1}$ and $\alpha \in \omega _{1}$, we have $G\cap D_{x}\neq \emptyset $ and $G\cap E_{x,\alpha }\neq \emptyset $. Let $F_{0}=\bigcup {G}$.  By Lemma  \ref{PRO_Lm18} there is a finitely independent, tight proper coloring $F$, such that $q\subset F_{0}\subset F$. The only thing left to prove is, that $p\subset F$. Choose any $x\in \mathscr{A}$. Then there is a unique $x'\in \mathscr{B}$ with $x{\equiv }^{*}x'$. Since $p$ is partially finitely independent, we have $p(x)=p(x')=q(x')=F(x')$. On the other hand $F$ is finitely independent, so $F(x)=F(x')$, thus $F(x)=p(x)$.
\end{proof}

\begin{corollary}\label{PRO_Cor6}
 If CH holds, then there is a finitely independent tight proper coloring $F:{}^{\omega }\omega _{1}\longrightarrow \omega _{1} $ that is not 2-tight.
\end{corollary}

\begin{proof}
 Let $z$ be the parity function, i.e. $z(i)=0$ if $i$ is even, $z(i)=1$, if $i$ is odd. Clearly for all $i\in \omega $, we have $z(i)\in \{ 0,1\} =\{ c_{0}(i),c_{1}(i) \}$. Let $\mathscr{A}=\{ c_{0},c_{1},z\} $, and $s$ be a function with $dom(s)=\mathscr{A}$, $s(c_{0})=0,s(c_{1})=1,s(z)=2$. Clearly $s$ has finite domain, $s$ is partially finitely independent and partially strongly proper. Then, by Theorem \ref{PRO_Thm8}, there is a finitely independent tight proper coloring with $F:{}^{\omega }\omega _{1} \longrightarrow \omega _{1}$ with $s\subset F$. We only need to show that $F$ is not 2-tight. If it was 2-tight, then by Lemma  \ref{PRO_Lm6} we would have $2=c(z)=F(z)\in \{F(c_{0}),F(c_{1})\} =\{s(c_{0}),s(c_{1})\} =\{ 0,1\}$ that is a contradiction.
\end{proof}

\section{Construction of tight proper colorings}

In this section we will construct some non-trivial tight proper colorings ${F:^{\omega}{\omega}\longrightarrow \omega }$. By Corollary  \ref{PRO_Cor5}, they cannot be finitely independent, and also they cannot be 2-tight.

\begin{definition}\label{PRO_Def10}
 A set $\mathscr{A}\subseteq {}^{\lambda }{\kappa }$ is \emph{lawful }if there are no $x,y\in \mathscr{A}$ with $x\perp y$.
\end{definition}

Clearly $F$ is a proper coloring if and only if $F^{-1}(\beta )$ is lawful for all $\beta \in C$.

\begin{lemma}\label{PRO_Lm19}
 An $F:{}^{\lambda }{\kappa } \longrightarrow C$ is a tight proper coloring if and only if $F^{-1}(\beta )$ is a maximal (by inclusion) lawful subset of ${}^{\lambda }{\kappa } $ for all $\beta \in C$.
\end{lemma}

\begin{proof}
 First suppose $F$ is a tight proper coloring. Then $F^{-1}(\beta )$ is a lawful subset of ${}^{\lambda }{\kappa } $ for all $\beta \in C$, and we need to show it is maximal. Suppose $F^{-1}(\beta )\subsetneq \mathscr{A}$ for some $\beta \in C$, where $\mathscr{A}$ is lawful and choose any $x\in \mathscr{A}\setminus F^{-1}(\beta )$. Then $F(x)\neq \beta $, so by tightness there is a $y\in {}^{\lambda }{\kappa } $, such that $y\perp x$ and $F(y)=\beta $. But then $y\in F^{-1}(\beta )\subseteq \mathscr{A}$ and also $x\in \mathscr{A}$, which contradicts that $\mathscr{A}$ is lawful. Thus, $F^{-1}(\beta )$ is maximal lawful for all $\beta \in C$.
 
 Now suppose $F$ has this property and choose any $x\in {}^{\lambda }{\kappa } $ and $\beta \in C$, with $F(x)\neq \beta $. Then $x\not\in F^{-1}(\beta )$ that is a maximal lawful set, so $F^{-1}(\beta )\cup \{ x\}$ is not lawful. Since $F^{-1}(\beta )$ is lawful, the only possible way is that there is some $y\in F^{-1}(\beta )$ that $y\perp x$. But then $F(y)=\beta $, so $y$ is a witness for the tightness in $x$ and $\beta $.
\end{proof}

\begin{lemma}\label{PRO_Lm26}
    If $T:{}^{\lambda }{\kappa }\rightarrow {}^{\lambda }{\kappa }$ is a bijection, such that for any $x,y\in {}^{\lambda }{\kappa }$, we have $T(x)\perp T(y)\Leftrightarrow x\perp y$, then for any subset $\mathscr{A}\subseteq {}^{\lambda }{\kappa }$, its image $T(\mathscr{A})$ is lawful if and only if $\mathcal{A}$ is lawful. The set $T(\mathscr{A})$ is maximal lawful iff $\mathscr{A}$ is maximal lawful.
\end{lemma}

\begin{proof}
    Straightforward.
\end{proof}

\begin{lemma}\label{PRO_Lm21}
 If $X\subseteq \lambda $ and $\mathscr{A}$ is a lawful subset of $^{X}\kappa $, then $${\mathscr{A}\times {}^{\lambda -X}\kappa =\{x\cup y:x\in \mathscr{A},y\in  {}^{\lambda -X}\kappa\}}$$ is a lawful subset of ${}^{\lambda }{\kappa } $. If $\mathscr{A}$ is maximal lawful, then $\mathscr{A}\times {}^{\lambda -X}\kappa $ is also maximal lawful.
\end{lemma}

\begin{proof}
 The first part of the lemma is straightforward.
 
 Now suppose $\mathscr{A}$ is maximal lawful. Then, by the first part, $\mathscr{A}\times {}^{\lambda -X}\kappa $ is lawful. To show it is maximal, suppose for contradiction that there is a lawful subset $\mathscr{B}$ of ${}^{\lambda }{\kappa } $, with $\mathscr{A}\times {}^{\lambda -X}\kappa \subsetneq \mathscr{B}$. Choose any $x\in \mathscr{B}\setminus \mathscr{A}\times {}^{\lambda -X}\kappa $. Then $x|_{X}\not\in \mathscr{A}$, and $\mathscr{A}$ is a maximal lawful subset of $^{X}\kappa $, so $\mathscr{A}\cup \{ x|_{X}\}$ is not lawful. But since $\mathscr{A}$ is lawful, the only possible way is that there is a $t\in \mathscr{A}$ with $t\perp x|_{X}$. Construct a $y\in {}^{\lambda }{\kappa } $ with $y|_{X}=t$ and $y\perp _{\lambda -X}x$. Then $y\perp x$ and $y\in \mathscr{A}\times {}^{\lambda -X}\kappa \subset \mathscr{B}$, but also $x\in \mathscr{B}$,which contradicts that $\mathscr{B}$ is lawful, so $\mathscr{A}\times {}^{\lambda -X}\kappa $ is maximal lawful.
\end{proof}

This way, we can extend tight proper colorings from $^{X}\kappa $ to ${}^{\lambda }{\kappa } $.

\begin{lemma}\label{PRO_Lm22}
 If $F:{}^{X}\kappa \longrightarrow C$ is a tight proper coloring, then the function $G:{}^{\lambda}\kappa \longrightarrow C$ defined by $G(x)=F(x|_{X})$ for $x\in {}^{\lambda }{\kappa } $ is also a tight proper coloring.
\end{lemma}

\begin{proof}
  For any $\beta \in C$, we have $G^{-1}(\beta )=\{ x\in {}^{\lambda }{\kappa } |G(x)=\beta \} =\{ x\in {}^{\lambda }{\kappa } |F(x|_{X})=\beta \} =\{ x\in {}^{\lambda }{\kappa } |x_{X}\in F^{-1}(\beta )\}=F^{-1}(\beta )\times {}^{\lambda -X}\kappa $. By Lemma  \ref{PRO_Lm19}. since $F$ is a tight proper coloring, $F^{-1}(\beta )$ is a maximal lawful subset of $^{X}\kappa $, and then by Lemma \ref{PRO_Lm21}. $G^{-1}(\beta )=F^{-1}(\beta )\times {}^{\lambda -X}\kappa $ is a maximal lawful subset of ${}^{\lambda }{\kappa } $ for all $\beta \in C$, thus by Lemma \ref{PRO_Lm19}. $G$ is a tight proper coloring.
\end{proof}

The simplest application of this lemma is when $X=\{ i\} $ for some $i\in \lambda $. Then a $\pi \in Sym(\kappa )$ clearly gives a function $F:{}^{\{ i\} }\kappa \longrightarrow \kappa $ with $F(x)=\pi (x) $. The extension of it by Lemma  \ref{PRO_Lm22} is the function $G:{}^{\lambda}\kappa \longrightarrow \kappa $ defined by $G(x)=\pi (x(i))$, the trivial coloring. In the next part we will construct some non-trivial tight proper colorings on $^{l}\omega $ and this way, we can extend them to get non-trivial tight proper colorings on $^{\omega }\omega $.

For any $x\in {}^{k}2$ let $\bar{x}\in {}^{k}2$ be such that $\bar{x}(i)=1-x(i)$ for all $i\in k$.

\begin{lemma}\label{PRO_Lm23}
 For any $k\in \omega $ and $\ell<2$ the set  
 \begin{displaymath}
    A_{k,\ell}=\Big\{ x\in {}^{2k+1}2\ |\ \sum_{i\in 2k+1}x(i)\equiv \ell(mod\:2)\Big\}
 \end{displaymath}
 is a  maximal lawful subset of $^{2k+1 }\omega $.

\end{lemma}

\begin{proof}
  It is enough to show that $A_{0}$ is maximal lawful, since the other case is similar.  For the lawfulness suppose for contradiction that $x\perp y$ for some $x,y\in A_{0}$. Since $x,y\in {}^{2k+1}2$, the only possible way would be that $y=\bar{x}$, but then $(\sum_{i\in 2k+1}x(i))+(\sum_{i\in 2k+1}y(i))=(\sum_{i\in 2k+1}x(i))+(\sum_{i\in 2k+1}\bar{x}(i))=(\sum_{i\in 2k+1}x(i))+(\sum_{i\in 2k+1}1-x(i))=\sum_{i\in 2k+1}x_{1}+(1-x(i))=2k+1$, so one of the sums is even, one is odd, $x,y\in A_{0}$ is not possible.
 
 For maximality suppose that $A_{0}\subsetneq B$, where $B$ is a lawful subset of $^{2k+1 }\omega $. Choose any $x\in B\setminus A_{0}$. If $x\in {}^{2k+1}2$, then $x\in A_{1}$, and as we have seen $\bar{x}\in A_{0}$. But then $x,\bar{x}\in B$ and $x\perp \bar{x}$, so $B$ is not lawful. Now we can suppose $x\not\in {}^{2k+1}2$, so there is an $i\in 2k+1$ that $x(i)\not\in 2$. We construct a $y\in A_{0}$ with $y\perp x$. For all $j\in 2k+1,j\neq i$, fix some $y(j)\in 2$ with $y(j)\neq x(j)$. There are 2 possible values for $y(i)$. Determine that value $y(i)\in 2$ to make $y\in A_{0}$. Since $y(i)\in 2$ and $x(i)\not\in 2$, we also have $y(i)\neq x(i)$, so $y\perp x$. But $x,y\in B$, so again $B$ is not lawful. Thus, $A_{0}$ is maximal lawful.
\end{proof}


\begin{theorem}\label{PRO_Thm9}
 For any $k\in \omega ,k\ge 1$, there is a non-trivial tight proper coloring $F:{}^{2k+1}\omega \longrightarrow \omega $ where all color classes are finite.
\end{theorem}

\begin{proof}
 The set $^{2k+1}\omega \times 2$ is countable, so we can fix a bijection $h:{}^{2k+1}\omega \times 2\longrightarrow \omega $. For $x\in {}^{2k+1}\omega$ let $F(x)=h(\lfloor \frac{x(0)}{2} \rfloor ,...,\lfloor \frac{x(2k)}{2}\rfloor ,\sum_{i\in 2k+1}x(i)(mod\: 2))$. We need to show that this is a tight proper coloring. To show this let $\beta \in \omega $ be arbitrary. Since $h$ is a bijection, there is a unique $l_{0},...,l_{2k}\in \omega $ and $\epsilon \in 2$, such that $h(l_{0},...,l_{2k},\epsilon )=\beta $. Then $F(x)=\beta $ for some $x\in {}^{2k+1}\omega $ if for all $i\in 2k+1$ we have $\lfloor \frac{x(i)}{2}\rfloor =l_{i}$, so $x(i)\in \{2l_{i},2l_{i}+1\}$ and also $\sum_{i\in 2k+1}x(i)=\epsilon $. Let $T_{l_{0},...,l_{2k}}:{}^{2k+1 }{\omega }\rightarrow {}^{2k+1 }{\omega }$, be the bijection defined by $(T(x))(i)=2l_{i}+x(i)$ for $i\in 2k+1$. The properties of Lemma \ref{PRO_Lm26} clearly hold for $T_{l_{0},...,l_{2k}}$.  Thus, we have $$F^{-1}(\beta )=\{x\in \bigtimes_{i\in 2k+1}\{ 2l_{i},2l_{i}+1\} |\sum_{i\in 2k+1}x(i)\equiv \epsilon (mod\: 2)\}=T_{l_{0},...,l_{2k}}(A_{k,\epsilon}).$$ By Lemma  \ref{PRO_Lm23} and lemma \ref{PRO_Lm26} this is a maximal lawful subset of $^{2k+1}\omega $, so by Lemma  \ref{PRO_Lm19}. $F$ is a tight proper coloring. The color classes are also finite, they have $2^{2k}$ elements.
\end{proof}



\begin{corollary}\label{PRO_Cor7}
 There are non-trivial tight proper colorings $F:{}^{\omega }\omega \longrightarrow \omega $.
\end{corollary}

\begin{proof}
 Let $k\in \omega ,k\ge 1$ arbitrary. By Theorem \ref{PRO_Thm9}, there is a  non-trivial tight proper coloring $G:{}^{2k+1}\omega \longrightarrow \omega $. Let $F(x)=G(x|_{2k+1})$. By Lemma \ref{PRO_Lm22}, $F$ is a tight proper coloring on $^{\omega }\omega $, and clearly $F$ is non-trivial, it depends on $2k+1$ coordinates.
\end{proof}

This way we can construct non-trivial tight proper colorings on $^{\omega }\omega $.

\begin{definition}\label{PRO_Def12}
A proper coloring $F:{}^{\lambda }{\kappa } \longrightarrow \mu $ {\em depends on finitely many coordinates} is there is an $S\in [\lambda ]^{<\omega }$ finite set, such that for all $x,y\in {}^{\lambda }{\kappa }$ with $x\equiv _{S}y$, we have $F(x)=F(y)$.
\end{definition}

All the colorings defined by theorem \ref{PRO_Thm9} depend on finitely many coordinates. Of course, we cannot make finitely independent $F$, but in the next part, we will construct some $F$, which does not depend on finitely many coordinates. For this, we will use coloring partitions.

\begin{definition}\label{PRO_C-tight}
A proper coloring $F:^{\lambda }{\kappa }\longrightarrow \mu $ is tight over $C\subseteq \mu $ if $Ran(F)\subseteq C$ and for all $x\in ^{\lambda }{\kappa }$ and $c\in C\setminus \{F(x)\}$, there is a $y\in ^{\lambda }{\kappa }$, with $y\perp x$ and $F(y)=c$.
\end{definition}

Clearly tight over $\mu $ simply means tight. It is also clear by Lemma $\ref{PRO_Lm18}$ that if ${F:^{\lambda }{\kappa }\longrightarrow \mu }$ is a proper coloring tight over $C$, then for all $c\in C$, the set $F^{-1}(c)\subseteq ^{\lambda }{\kappa }$ is maximal lawful. 

\begin{lemma}\label{PRO_Lm20}
 If $F:{}^{\lambda }{\kappa } \longrightarrow C_{1}$ is a proper coloring tight over $C_{1}$ and $h:C_{1}\longrightarrow C_{2}$ is a bijection, then $h\circ F:{}^{\lambda}\kappa \longrightarrow C_{2}$ is a proper coloring tight over $C_{2}$.
\end{lemma}

\begin{proof}
Straightforward.
\end{proof}

\begin{definition}\label{PRO_Def11}
  $\mathscr{P}$ is a \emph{coloring partition on ${}^{\lambda }{\kappa } $ with index set $I$} if $\mathscr{P}=\{ (\mathscr{A}_{i},F_{i},C_{i},B_{i}):i\in I\} $ is a set of indexed quadruples, where $\mathscr{A}_{i}\subseteq {}^{\lambda }{\kappa } $, $F_{i}:{}^{\lambda}\kappa \longrightarrow C_{i}$ is a proper coloring and $B_{i}\subseteq C_{i}$ for all $i\in I$, and the following properties hold:
\begin{enumerate}[1.] 
\item   $\mathscr{A}_{i}\cap \mathscr{A}_{j}=\emptyset $ for $i,j\in I,i\neq j$,
\item
  $\bigcup_{i\in I}\mathscr{A}_{i}={}^{\lambda }{\kappa } $,
 \item
  $F_{i}$ does not split $\mathscr{A}_{i}$ for any $i\in I$ 
  (i.e. for all $\beta \in C_{i}$ either $F_{i}^{-1}(\beta )\subseteq \mathscr{A}_{i}$ or $F_{i}^{-1}(\beta )\cap \mathscr{A}_{i}=\emptyset $), 
 \item 
  $Ran(F_{i}|_{\mathscr{A}_{i}})=B_{i}$ for $i\in I$,
\item 
  $B_{i}\cap B_{j}=\emptyset $ for $i,j\in I,i\neq j$.
\end{enumerate}

 By properties 1 and 2 the family $\{\mathscr{A}_{i}:i\in I\}$ is a partition of ${}^{\lambda }{\kappa } $. We define the \emph{induced function} $F$ with $F|_{\mathscr{A}_{i}}=F_{i}|_{\mathscr{A}_{i}}$, this is a well-defined function.
\end{definition}

Clearly for a coloring partition $\mathscr{P}$, the range of the induced function is $Ran(F)=\bigcup_{i\in I}Ran(F|_{\mathscr{A}_{i}})=\bigcup_{i\in I}Ran(F_{i}|_{\mathscr{A}_{i}})=\bigcup_{i\in I}B_{i}$.

\begin{lemma}\label{PRO_Lm24}
 If $\mathscr{A}$ is a coloring partition on ${}^{\lambda }{\kappa } $ with index set $I$, then the induced function $F$ is a proper coloring. If $F_{i}:{}^{\lambda}\kappa \longrightarrow C_{i}$ is a $C_{i}$-tight proper coloring for all $i\in I$, $B=\bigcup_{i\in I}B_{i}$, then $F:{}^{\lambda}\kappa \longrightarrow B$ is $B$-tight.
\end{lemma}

\begin{proof}
 We prove this by examination of color classes. Let $\beta \in B=\bigcup_{i\in I}B_{i}$. Then there is some $i\in I$, such that $\beta \in B_{i}$. First we will show that $F^{-1}(\beta )\subseteq \mathscr{A}_{i}$. Suppose for contradiction that there is some $x\in F^{-1}(\beta )$ with $x\not\in \mathscr{A}_{i}$. Then $x\in \mathscr{A}_{j}$ for some $j\in I,j\neq i$, thus $F(x)=F_{j}(x)\in Ran(F_{j}|_{\mathscr{A}_{j}})=B_{j}$. But $F(x)=\beta \in B_{i}$ and by property 5 $B_{i}\cap B_{j}=\emptyset $, so that is a contradiction. Then we have $F^{-1}(\beta )=\{ x\in \mathscr{A}_{i}| F(x)=\beta \} =\{ x\in \mathscr{A}_{i}| F_{i}(x)=\beta \} =\mathscr{A}_{i}\cap F_{i}^{-1}(\beta )$. By property 3 either $F_{i}^{-1}(\beta )\subseteq \mathscr{A}_{i}$ or $F_{i}^{-1}(\beta )\cap \mathscr{A}_{i}=\emptyset $, but the latter case is not possible, since $\beta \in B_{i}=Ran(F_{i}|_{\mathscr{A}_{i}})$, so the former case holds. Thus, $F^{-1}(\beta )=F_{i}^{-1}(\beta )$. Since $F_{i}$ is a proper coloring and $\beta \in B_{i}\subseteq C_{i}$ this set is lawful for all $\beta \in B$, and if the functions $F_{i}$ are tight, then by Lemma  \ref{PRO_Lm19}. they are maximal lawful. Thus, $F$ is a proper coloring and if all $F_{i}$-s are tight, by Lemma \ref{PRO_Lm19}. $F$ is also tight. 
\end{proof}

\begin{theorem}\label{PRO_Thm10}
 There is a tight proper coloring $F:{}^{\omega }\omega \longrightarrow \omega $ that does not depend on just finitely many coordinates.
\end{theorem}

By Corollary \ref{PRO_Cor5}, $F$ cannot be finitely independent.

\begin{proof}
 We construct a coloring partition $\mathscr{P}$ on $^{\omega }\omega $ with the index set $I=\omega $. For all $i\in \omega $ let $\mathscr{A}_{i}=\{ x\in {}^{\omega }\omega |x(0)\in \{ 2i,2i+1\} \}$. Let $(C_{i})_{i\in \omega }$ be a partition of $\omega $ into infinite pieces.
 For all $i\in \omega $, let $h_{i}:{}^{2i+1}\omega \times 2\longrightarrow C_{i}$ be a bijection and let us define the functions $G_{i}:^{2i+1}{\omega}\longrightarrow C_{i}$ by $G_{i}(x)=h_{i}(\lfloor \frac{x(0)}{2} \rfloor ,...,\lfloor \frac{x(2i)}{2}\rfloor ,\sum_{j\in 2i+1}x(j)(mod\: 2))$. By Lemma \ref{PRO_Lm23} that for all $c\in C_{i}$, the set $G_{i}^{-1}(c)$ is maximal lawful, so by Lemma \ref{PRO_Lm19} $G_{i}$ is $C_{i}$-tight.
 By Lemma \ref{PRO_Lm22}. we can extend them to $C_{i}$ tight proper colorings $F_{i}:{}^{\omega }\omega \longrightarrow C_{i}$ with $F_{i}(x)=G_{i}(x|_{2i+1})$ for all $x\in {}^{\omega }\omega $. Also let $B_{i}=Ran(F_{i}|_{\mathscr{A}_{i}})\subseteq Ran(F_{i})=C_{i}$.
 
 Now let us prove that $\mathscr{P}=\{ (\mathscr{A}_{i},F_{i},C_{i},B_{i}):i\in \omega \}$ is a coloring partition. Properties 1. and 2. are trivial, since the sets $\mathscr{A}_{i}$ form a partition depends only on the first coordinate. For property 3 choose an $x\in \mathscr{A}_{i}$ and $y\not\in \mathscr{A}_{i}$, we need that $F_{i}(x)\neq F_{i}(y)$. Clearly $\lfloor \frac{x(0)}{2}\rfloor =i$ and $\lfloor \frac{y(0)}{2}\rfloor \neq i$, so $\lfloor \frac{x(0)}{2}\rfloor \neq \lfloor \frac{y(0)}{2}\rfloor $. Then $F_{i}(x)=G_{i}(x|_{2i+1})=h_{i}(\lfloor \frac{x(0)}{2}\rfloor ,...,\lfloor \frac{x(2i)}{2}\rfloor ,\sum_{j\in 2i+1}x(j)(mod\: 2))\neq h_{i}(\lfloor \frac{y(0)}{2}\rfloor ,...,\lfloor \frac{y(2i)}{2}\rfloor ,\sum_{j\in 2i+1}y(j)(mod\: 2))=G_{i}(y|_{2i+1})=F_{i}(y)$. 
 Property 4 is trivial by the definition of $B_{i}$. For property 5 if $i\neq j$, then $B_{i}\cap B_{j}\subseteq C_{i}\cap C_{j}=\emptyset $. Thus, all properties hold, $\mathscr{P}$ is a coloring partition.
 
 Let $G:{}^{\omega }\omega \longrightarrow B$ be the function induced by $\mathscr{P}$, where $B=\bigcup_{i\in \omega }B_{i}$. By Lemma \ref{PRO_Lm24}. this is a $B$-tight proper coloring. Clearly $B\subseteq \omega $, on the other hand $B_{i}\neq \emptyset $ for all $i\in \omega $, and they are pairwise disjoint, so $B$ is countably infinite. Choose a bijection $h:B\longrightarrow \omega $ and let $F=h\circ G$. By Lemma  \ref{PRO_Lm20}. $F: {}^{\omega }\omega \longrightarrow \omega $ is a tight proper coloring. The only thing that is left for proof that $F$ does not depend only on a finite set of coordinates.
 
 For this let $S\subset \omega $ be an arbitrary finite set and we need to construct some $x,y\in {}^{\omega }\omega $ with $x\equiv _{S}y$, but $F(x)\neq F(y)$. Choose any $i\ge 1$ with $2i>max(S)$. Then clearly $2i\not\in S$. Let $x(0)=y(0)=2i$, $x(2i)=0$, $y(2i)=1$ and $x(j)=y(j)=0$ for $j\not\in \{ 0,2i\}$. Then $x,y$ differs only in the coordinate $2i$, so $x\equiv _{S}y$ clearly holds. Since $\sum_{j\in 2i+1}x(j)=2i\equiv 0(mod\: 2)$ and $\sum_{j\in 2i+1}y(j)=2i+1\equiv 1(mod\: 2)$, we have $F_{i}(x)=G_{i}(x|_{2i+1})=h_{i}(i,0,...,0,0)\neq h_{i}(i,0,...,0,1)=G_{i}(y|_{2i+1})=F_{i}(y)$. Also since $x,y\in \mathscr{A}_{i}$, we have $G(x)=F_{i}(x)\neq F_{i}(y)=G(y)$, and since $h$ is a bijection, we have $F(x)=h(G(x))\neq h(G(y))=F(y)$. Thus, $F$ does not only depend on any finite set $S$.
\end{proof}

By similar construction we can make tight proper colorings on $^{\omega }\omega _{1}$ and in general on any ${}^{\lambda }{\kappa } $ that are neither trivial nor depend only on a finite set of coordinates, but they are not finitely independent. As 2-tight proper colorings are either finitely independent (if the corresponding ultrafilter is non-principal) or trivial (if the corresponding ultrafilter is principal), this construction gives tight proper colorings that are not 2-tight.

\bibliographystyle{plain}

\end{document}